\documentclass[a4paper,12pt]{article}%
\usepackage{CJK}
\usepackage{setspace}
\usepackage{amsfonts}
\usepackage{makeidx}
\usepackage{latexsym}
\usepackage{amsmath}
\usepackage{color}
\usepackage{array}
\usepackage{booktabs}
\usepackage{multirow}
\usepackage{amssymb}
\usepackage{graphicx, amsmath, amsthm, amssymb, bm, boxedminipage}
\usepackage{geometry}
\usepackage{graphicx}%
\allowdisplaybreaks[4]
\allowdisplaybreaks[4]
\newtheorem{example}{Example}
\numberwithin{equation}{section}
\theoremstyle{plain}
\newtheorem{thm}{Theorem}[section]
\newtheorem{lemma}[thm]{Lemma}
\newtheorem{remark}[thm]{Remark}

\newtheorem{corol}[thm]{Corollary}

\ifx\pdfoutput\relax\let\pdfoutput=\undefined\fi
\newcount\msipdfoutput
\ifx\pdfoutput\undefined\else
\ifcase\pdfoutput\else
\ifx\paperwidth\undefined\else
\ifdim\paperheight=0pt\relax\else\pdfpageheight\paperheight\fi
\ifdim\paperwidth=0pt\relax\else\pdfpagewidth\paperwidth\fi
\fi\fi\fi
\begin{document}

\title{A Central Limit Theorem, Loss Aversion and Multi-Armed Bandits\thanks{Chen is
at School of Mathematics, Shandong University, zjchen@sdu.edu.cn, Epstein is
at Department of Economics, McGill University, larry.epstein@mcgill.ca, and
Zhang is at School of Mathematics, Shandong University, zhang$\_$%
gd@mail.sdu.edu.cn. Chen gratefully acknowledges the support of the National
Key R\&D Program of China (grant No. ZR2019ZD41), and the Taishan Scholars
Project. We are grateful to Peter Wakker for helpful discussions and
references.}}
\author{Zengjing Chen
\and Larry G. Epstein
\and Guodong Zhang}
\date{May 15, 2022}
\maketitle

\begin{abstract}
This paper studies a multi-armed bandit problem where the decision-maker is
loss averse, in particular she is risk averse in the domain of gains and risk
loving in the domain of losses. The focus is on large horizons. Consequences
of loss aversion for asymptotic (large horizon) properties are derived in a
number of analytical results. The analysis is based on a new central limit
theorem for a set of measures under which conditional variances can vary in a
largely unstructured history-dependent way subject only to the restriction
that they lie in a fixed interval.

\bigskip\vspace{0.4in}


\end{abstract}

\newpage

\section{Introduction}

We study the following (multi-armed bandit) sequential choice
problem.\footnote{Overviews and textbook treatments of the bandit model can be
found in Berry and Fristadt (1985) and Slivkins (2019), for example. The first
application to economics is Rothschild (1974). See Bergemann and Valimaki
(2008) for references to a range of other economic applications.} There are
finitely many arms (or actions), each yielding a random payoff. Probability
distributions have a common mean but differ otherwise and may not be known to
the decision-maker (DM). At each stage $i=1,2,...,n$, DM chooses one arm,
knowing the realized outcomes from previous choices. Ex ante she chooses a
strategy to maximize expected utility, where the utility index is a function
of the (suitably weighted) average payoff. Because we are interested in
varying horizons, it is convenient to define a strategy for an infinite
horizon, and then to use its truncation for any given finite horizon. Refer to
a strategy as \emph{asymptotically optimal }if the expected utility it implies
in the limit as horizon $n\rightarrow\infty$ is at least as large as that
implied by any other strategy. We study large-horizon approximations to the
value (indirect utility) of the bandit problem and corresponding
asymptotically optimal strategies.

A second novelty in our model is the assumption that DM is loss averse (global
risk aversion is a limiting special case). Loss aversion was introduced via
cumulative prospect theory by Tversky and Kahneman (1992), and has since been
well-established empirically and widely applied in economics and finance (see
for example, Kahneman and Tversky 2000, Kobberling and Wakker 2005, Barberis
2013, and the references therein). Its essential elements are (i) a reference
point; (ii) utility depends only on gains and losses relative to that
reference point rather than on the total payoff (or total wealth); (iii) risk
aversion (concavity) for gains and risk loving (convexity) for losses; and
(iv) greater sensitivity to losses than to gains. Our interest in this paper
is the effect of loss aversion in the sequential context defined by a bandit
problem. To our knowledge, this is the first study of loss aversion in bandit
problems.\footnote{Xu and Zhou (2013) and Ebert and Strack (2015) study
optimal stopping problems assuming prospect theory. Their focus is on the
probability-weighting aspect of prospect theory and loss aversion plays no
role in their analyses. Two studies of loss aversion in a sequential context
are Easley and Yang (2015) and Shi et al (2015). The former uses numerical
analysis to study the wealth and price effects of loss aversion in the
equilibrium of a dynamic heterogeneous-agent economy. The latter models
dynamic portfolio choice with loss aversion, where the reference point varies
endogenously in response to prior wealth outcomes. In both cases, analysis is
largely numerical and there is little overlap with the bandit literature in
general, and with our model, in particular. Guasoni et al (2020) study
shortfall aversion, which shares the spirit of loss aversion but which is more
directly relevant to preference over deterministic consumption streams rather
than over lotteries.}

We have two related reasons for studying asymptotics. First, it promotes
tractability and the derivation of analytical results. Though the literature
on bandit problems is enormous, theoretical analysis of Bayesian models is, to
the best of our knowledge, restricted to the special case of risk neutrality
(see section \ref{section-optimization} for elaboration and a
qualification).\footnote{Two studies of bandit problems that explicitly
address risk are Sani, Lazaric and Munos (2013) and Huo and Fu (2017). They
assume regret minimization rather than expected utility maximization, and
focus on computational algorithms rather than on qualitative theoretical
results.} Besides its obvious limitations, risk neutrality also imposes the
invariance of risk attitude as some outcomes are realized, and this invariance
is key to well-known sequential properties of optimal strategies derived in
the literature.\footnote{For example, in an infinite-horizon setting where
means can differ, and with one unknown arm and one arm whose distribution is
known, then once the known arm is chosen it will continue to be chosen
therafter (Rothschild 1974, pp. 190-191).} In contrast, endogenously varying
risk attitude is inherent in loss aversion. Moreover, in our setting where
means are known and common to all arms, risk neutrality would trivialize the problem.

Our second reason for studying asymptotics is that tractability is plausibly a
concern not only for the modeler but also for the decision-maker within the
model. We view her as struggling to comprehend an extremely complicated
finite-horizon optimization problem, and adopting instead the simplifying
assumption of an infinite horizon. She does so with the recognition that an
asymptotically optimal strategy is approximately optimal if her horizon is
sufficiently long.

Here is an informal outline of some of our analytical results, which obtain as
stated in the infinite-horizon limit and approximately for sufficiently large
finite horizons.

\begin{enumerate}
\item Maximum ex ante expected utility depends on the distributions describing
each arm \emph{only through their variances}. Moreover, it depends only on the
largest and smallest variance. Consequently, at each history, arms with
intermediate conditional variances can be ignored.

\item Depending on the reference point, it is possible to achieve a level of
ex ante expected utility that is equal to, or strictly greater than, the level
when the payoff to each arm is riskless. In that sense, \emph{risk may be
desirable} in the sequential context, even though "comparable" risks would be
rejected in a one-shot choice setting.

\item Suppose that the distributions describing every arm are known. Then, in
spite of the absence of learning, an asymptotically optimal strategy switches
indefinitely between two fixed extreme arms (those with the smallest\ and
largest variances) as the decision-maker moves between cumulative gains and
cumulative losses. Given two arms that exhibit the two extreme variances, all
other arms are redundant.

\item Suppose there are two arms and that the pair of variances is known, but
there is prior uncertainty about which arm has which variance. Then it is
asymptotically optimal to choose myopically at each stage, that is, as though
there are no subsequent choices to be made.

\item None of the above results rely on assumptions about the nature of risk
aversion in the domain of gains or about the nature of risk loving in the
domain of losses. They depend only on preference over "mixed" lotteries.
\end{enumerate}

Finally, we turn to the proofs of these and other results about bandits and
loss aversion. It is not surprising that asymptotic results may be approached
via limit theorems. However, classic limit results do not apply, and the key
to our proofs is a new central limit theorem (CLT). The martingale version of
the central limit theorem considers a sequence $\left(  X_{i}\right)  $ of
random variables having zero conditional mean and constant conditional
variance $\sigma^{2}$, and shows that (under suitable additional conditions)
the distribution of $\Sigma_{i=1}^{n}X_{i}/\sqrt{n}$ converges to the normal
$\mathbb{N}\left(  0,\sigma^{2}\right)  $ as $n\rightarrow\infty$. (The
classic result for identically and independently distributed random variables
is an immediate special case). This paper establishes a CLT under the relaxed
assumption on variance according to which conditional variances can vary in a
largely unstructured history-dependent way subject only to the restriction
that they lie in a fixed interval $\left[  \underline{\sigma}^{2}%
,\overline{\sigma}^{2}\right]  $, in which case limits take a novel and
tractable form. This CLT is the main technical contribution of the paper. One
well-known motivation for generalizing from a single probability distribution
(hence single variance) to a set of probability distributions (hence set of
variances) is robustness to model uncertainty or ambiguity. However, model
uncertainty plays no role in our bandit problem - DM is a Bayesian agent,
perfectly confident in her understanding of the environment - thus
highlighting the usefulness of sets of measures even for Bayesian models.

We proceed as follows. The bandit model and the results outlined above are
described in detail in the next section. Proofs for these results must await
the CLT which is presented next in section \ref{section-CLT}. Proofs of the
CLT and related results are presented in Appendix \ref{app-CLT} and proofs for
the bandit application are in Appendix \ref{app-bandits}.

\section{Multi-Armed Bandits\label{section-bandits}}

\subsection{Beliefs, utility and optimization}

Let $\mathcal{A}$ be a finite set of arms (or actions). The outcome of any
action lies in the finite set $\overline{\Omega}\subset\mathbb{R}$. Thus
outcome sequences lie in $\Omega=\Pi_{1}^{\infty}\Omega_{i}$, where
$\Omega_{i}=\overline{\Omega}$ for each $i$. The timing is as follows: At each
$i\geq1$, the history $\omega^{\left(  i-1\right)  }=\left(  \omega
_{1},...,\omega_{i-1}\right)  $ is known\textit{, (}$\omega^{\left(  0\right)
}=\varnothing$), an action $a_{i}\in\mathcal{A}$\textit{ }is chosen, and then
the resulting outcome $\omega_{i}$ is realized. Define $X_{i}\left(
\omega\right)  =\omega_{i}$, the outcome at stage $i$.

Let $\mathcal{G}_{i-1}$ be the $\sigma$-algebra representing information at
stage $i$, ($\mathcal{G}_{0}=\{\Omega,\varnothing\}$), and let $\mathcal{G}%
=\sigma(\cup_{1}^{\infty}\mathcal{G}_{i})$ be the corresponding $\sigma
$-algebra on $\Omega$.

The outcome resulting from any action is uncertain and the choice of a
contingent plan, or strategy, is determined by expected utility maximization.
The remaining primitives of the model - strategies, beliefs and the vNM
utility index - are described next.

\subsubsection{Strategies and beliefs\label{section-beliefs}}

The contingent choice of action at stage $i$ depends on (conditional) beliefs
about the next outcome, which generally depend on the arm being considered and
also on what is learned from previous choices and their outcomes. Importantly,
the inference to be drawn from the history $\omega^{\left(  i-1\right)  }$ of
outcomes depends on which arms produced them. Thus, the choice of action at
stage $i$ is expressed as
\begin{equation}
a_{i}={s_{i}(a^{(i-1)},\omega^{(i-1)})}\text{,} \label{ai}%
\end{equation}
where $a^{\left(  i-1\right)  }=\left(  a_{1},...,a_{i-1}\right)  $ denotes
the history of past actions ($a^{0}=\varnothing$). Refer to $s_{i}%
:\mathcal{A}^{i-1}\times\Pi_{j=1}^{i-1}\Omega_{j}\rightarrow\mathcal{A}$ as
{the \emph{strategy at stage }$i$, and denote the set of all such }$s_{i}$ by
$\mathcal{S}_{i}$.{ The infinite sequence $s=\left(  s_{i}\right)
_{1}^{\infty}$ is called simply a \emph{strategy}. }The corresponding set of
strategies is $\mathcal{S}$.\footnote{For any given $n$, $s\in\mathcal{S}$
induces the contingent plan $\left(  s_{i}\right)  _{1}^{n}$, which is
adequate if one is interested only in the $n$-horizon case. Because we will be
interested in varying horizons, it is convenient to define a strategy to apply
to all finite horizons.}

Turn to beliefs. For the reasons noted above, beliefs about the next outcome
depend on both the action being considered, hence on the strategy for the
current stage, and on the history of past actions. Thus we model these beliefs
for stage $i$ by the conditional probability measure
\begin{equation}
P_{i}^{s_{i}}=P_{i}^{s_{i}}\left(  \cdot\mid a^{(i-1)},\omega^{(i-1)}\right)
\in\Delta\left(  \Omega_{i}\right)  \text{.} \label{Psi}%
\end{equation}
The set of $1$-step-ahead conditionals $\left\{  P_{i}^{s_{i}}\right\}
_{i\geq1,~s_{i}\in\mathcal{S}_{i}}$ is a \emph{primitive} that represents
beliefs (which may be taken to be subjective or objective).

Given a (fixed) strategy {$s=(s_{1},...,s_{i},...)$, we can combine the
primitive conditionals into a measure }$P^{s}$ on $\Pi_{1}^{\infty}\Omega_{i}%
$. To describe why and how, suppose that DM is considering the strategy $s$.
Then she knows that the action history $a^{(i-1)}$ at any future stage $i$ is
uniquely determined by the corresponding outcome history $\omega^{\left(
i-1\right)  }$, and by the given stage strategies $s_{1},...,s_{i-1}$, via
repeated iteration of the relation%
\begin{equation}
a_{j}=s_{j}\left(  a^{\left(  j-1\right)  },\omega^{\left(  j-1\right)
}\right)  \text{, }j=1,2,...,i-1\text{\textbf{.}} \label{asj}%
\end{equation}
Consequently, DM faces uncertainty only about future outcomes and she seeks a
probability measure over all outcome sequences that is consistent with the
primitive 1-step-ahead conditionals $\left\{  P_{i}^{s_{i}}\right\}
_{i\geq1,~s_{i}\in\mathcal{S}_{i}}$. As noted, the action history appearing in
(\ref{Psi}) can be substituted out and the conditioning information becomes a
function of $\omega^{\left(  i-1\right)  }$ alone. Therefore, the 1-step-ahead
conditionals can be pasted together in the usual fashion: By the
Ionescu-Tulcea extension theorem, one obtains a (unique) measure $P^{s}$,
\begin{equation}
P^{s}\in\Delta\left(  \Pi_{1}^{\infty}\Omega_{i},\mathcal{G}\right)  \text{.}
\label{Ps}%
\end{equation}
Moreover, its 1-step-ahead conditional $P_{i}^{s}\left(  \cdot\mid
\mathcal{G}_{i}\right)  \in\Delta\left(  \Omega_{i}\right)  $
\textquotedblleft agrees" with the primitive conditional $P_{i}^{s_{i}}$ in
the sense that%
\begin{equation}
P_{i}^{s}\left(  \cdot\mid\mathcal{G}_{i}\right)  \left(  \omega^{\left(
i-1\right)  }\right)  =P_{i}^{s_{i}}\left(  \cdot\mid a^{(i-1)},\omega
^{(i-1)}\right)  \text{,} \label{agree}%
\end{equation}
where $a^{\left(  i-1\right)  }$ is obtained from (\ref{asj}).

We assume that each $P_{i}^{s}$ has \emph{full support} on $\overline{\Omega}%
$. We assume also that mean outcomes are common to all arms (hence also
strategies) and fixed:
\begin{equation}
E_{P^{s}}\left[  X_{i}\mid\mathcal{G}_{i-1}\right]  =m=0\text{ \ for all
}i\geq1\text{ and all }s\in\mathcal{S}\text{,} \label{muzero1}%
\end{equation}
where setting $m=0$ is without loss of generality. Finally, the extreme
(largest and smallest) variances across all arms are assumed constant:
\begin{equation}
ess\sup\limits_{s\in\mathcal{S}}E_{P^{s}}\left[  X_{i}^{2}|\mathcal{G}%
_{i-1}\right]  =\overline{\sigma}^{2}\text{ and }ess\inf_{s\in\mathcal{S}%
}E_{P^{s}}\left[  X_{i}^{2}|\mathcal{G}_{i-1}\right]  =\underline{\sigma}%
^{2}\text{ for all }i\geq1\text{,} \label{condvar1}%
\end{equation}
for some $\overline{\sigma}^{2}>\underline{\sigma}^{2}>0$.

Apart from these conditions, the history-dependence of any primitive
conditional $P_{i}^{s_{i}}$ is unrestricted, and consequently so is the
implied model of learning. For example, the model does not rule out that past
outcomes of arm $a$ may inform beliefs about future outcomes of arm
$a^{\prime}\not =a$. Independence across arms, or other hypotheses about the
nature of learning, can be accommodated by suitably specializing beliefs, but
are not assumed in the general model. (Sections \ref{section-nolearning} and
\ref{section-2arms} describe two such specializations.)

For readers who find the strategy-dependence of probability measures
unorthodox we add that it is readily understood in the following terms.
Consider a generic static choice problem of the form $\sup_{a\in\mathcal{A}%
}E_{\mu}\left[  u\left(  X^{a}\right)  \right]  $, where $X^{a}$ is the random
variable outcome associated with action $a$ and $\mu$ is a prior over the
underlying state space $\Omega$.\footnote{Here $\Omega$ is an abstract state
space, not necessarily related to the product state space used in the bandit
model. Similarly, for $\mathcal{A}$ and for $X$ below.} Then each $X^{a}$
induces a probability distribution, denoted $p^{a}$, over $\Omega$, and the
preceding optimization problem can be written as $\sup_{a\in\mathcal{A}%
}E_{p^{a}}\left[  u\left(  X\right)  \right]  $, where $X\left(
\omega\right)  =\omega$. Thus the choice between actions, modeled as the
choice between random variables, can be expressed alternatively as the choice
between action-dependent probability distributions over outcomes (that is,
lotteries). The analogue of this reformulation for our sequential choice
context leads to strategy-dependent probability measures.\footnote{The use of
action-dependent probabilities (or moral hazard) has been recognized in the
decision theory literature (Dreze 1987, Kelsey and Milne 1999, and Karni 2011,
for example). These papers are concerned primarily with axiomatic foundations,
extending those for subjective expected utility, while our motivation in
studying the bandit problem is more applied. We differ also in our focus on
sequential choice.}

\subsubsection{Utility\label{section-utility}}

We assume that, at each stage $i$, outcomes for each action are evaluated
according to whether they produce gains or losses relative to a reference
point, which we take to be their common mean (taken to be zero for
convenience). Then $X_{i}$ gives the gain/loss at stage $i$. Since
gains/losses are incurred at each stage, they must be aggregated. We posit
that, for any horizon length $n$, utility depends on their $\sqrt{n}$-weighted
average. Consequently, given the strategy $s$, the implied stream of
gains/losses has expected utility given by%
\begin{equation}
U_{n}\left(  s\right)  =E_{P^{s}}\left[  \varphi\left(  \Sigma_{1}^{n}%
X_{i}/\sqrt{n}\right)  \right]  \text{,} \label{Un}%
\end{equation}
where $\varphi$ is the vNM utility index, which will be described shortly.

The $\sqrt{n}$-weighted averaging calls for some discussion. Consider a
setting (such as a casino, where trials correspond to playing one or another
slot machine or gambling device) where the time between trials is so small as
to preclude discounting, and where the monetary payoffs at different trials
are perfect substitutes. We are not aware of any axiomatic (or empirical)
guidance for how a decision-maker does or should aggregate or average money
streams in this context given that arbitrarily large horizons are relevant.
The unweighted arithmetic average might be slightly simpler to contemplate and
calculate, but significantly, it also reflects a specific and possibly
inappropriate weighting to finite sets of trials. Indeed, as is familiar from
discussions of the classic law of large numbers (LLN) and CLT, one might argue
that scaling by $\frac{1}{n}$ implies "too little" weight for finite sets of
trials, particularly when considering volatility. That will be argued also in
our setting (Remark \ref{remark-LLN}) and is our practical modelers' rationale
for preferring (\ref{Un}).

\begin{remark}
To be perfectly clear, the utility functions $U_{n}$ rank strategies for any
given horizon $n$. They do not rank horizons. That is, statements such as
$U_{n}\left(  s\right)  \geq U_{n}\left(  s^{\prime}\right)  $ are meaningful,
but statements such as $U_{n}\left(  s\right)  \geq U_{n^{\prime}}\left(
s\right)  $ are not and do not play a role below.
\end{remark}

The utility index $\varphi$ appearing in (\ref{Un}) is defined by
\begin{equation}
\varphi(x)=\left\{
\begin{array}
[c]{lll}%
\ \ \varphi_{1}(x-c) &  & x\geq c\\
-\theta^{-1}\varphi_{1}\left(  -\theta(x-c)\right)  &  & x<c
\end{array}
\right.  \label{phi-trinomial}%
\end{equation}
where we assume:

\smallskip

\noindent\textbf{Assumption-Utility}: $\theta=\underline{\sigma}%
/\overline{\sigma}<1$, $\varphi_{1}(0)=0$, $\varphi_{1}\in C_{b}^{3}\left(
\mathbb{R}_{+}\right)  $, and $\varphi_{1}$ is (strictly) increasing and
(strictly) concave for $x>c$.\footnote{$C_{b}^{3}\left(  \mathbb{R}%
_{+}\right)  $ is the set of functions on the non-negative real line with
continuous and bounded third order derivatives.}

\smallskip

\noindent\noindent Then, $\varphi$ is increasing globally, concave for $x>c$
(corresponding to gains) and convex for $x<c$ (corresponding to losses),
implying risk aversion for gains and risk seeking for losses. In addition,
\begin{equation}
x>y\geq0\text{ ~}\Longrightarrow\text{ ~}\left(  c+y,\tfrac{1}{2}%
;c-y,\tfrac{1}{2}\right)  \succ\left(  c+x,\tfrac{1}{2};c-x,\tfrac{1}%
{2}\right)  \text{,} \label{LAdefn}%
\end{equation}
indicating greater sensitivity to the increased loss ($-x<-y$) than to the
increased gain ($x>y$). In differential form, it states that
\begin{equation}
\varphi^{\prime}\left(  c-x\right)  >\varphi^{\prime}\left(  c+x\right)
\text{, \ for all }x>0\text{.} \label{LA}%
\end{equation}
We take these to be the defining properties of (strict) loss aversion,
following Wakker and Tversky (1993, p. 164), for example. An implication is
that $-\varphi\left(  c-x\right)  >\varphi\left(  c+x\right)  $, \ for all
$x>0$, that is, the lottery $\left(  c+x,\frac{1}{2};c-x,\frac{1}{2}\right)  $
is strictly inferior to receiving $0$ for sure.

The following example will be useful in the sequel (see (\ref{Vexp})) because
of its tractability.

\begin{example}
[Exponential]\label{example-exp}Let $\varphi_{1}\left(  x\right)
=1-\exp\left(  -x\right)  $, so that
\begin{equation}
\varphi(x)=\left\{
\begin{array}
[c]{ll}%
1-\exp(-(x-c)) & x\geq c\\
\theta^{-1}(\exp(\theta(x-c))-1) & x<c
\end{array}
\right.  \label{exp}%
\end{equation}
where $c\in\mathbb{R}$ and $\theta=$\underline{$\sigma$}/$\overline{\sigma}$.
\end{example}

Because of its origins in prospect theory, loss aversion is often viewed as
tied to probability weighting or distortion, (which is absent in our expected
utility model), and also to a kink in the utility index at the reference point
(which is also absent here because $\varphi$ defined above is continuously
differentiable everywhere). However, neither is necessary mathematically or
conceptually for the above behavioral properties that define loss
aversion.\footnote{Kobberling and Wakker (2005) argue explicitly for a
conceptual separation between loss aversion and probability weighting. They
write (p. 124): \textquotedblleft We have introduced utility, probability
weighting and loss aversion as logically independent factors of risk attitude
... their (in)dependence empirically is more intricate."} Accordingly,
consistent with common practice, we exclude probability distortions, hence
Allais-type behavior, in order to isolate the effects of loss aversion on
sequential decision-making. As for a kink, it has limited empirical content;
for example, a finite set of pairwise rankings of lotteries, as is common in
experimental investigations of loss aversion, cannot refute differentiability.
Moreover, the theoretical connection of a kink to loss aversion is very much
dependent on the choice of functional form. For example, suppose that, instead
of (\ref{phi-trinomial}), one posits that
\begin{equation}
\varphi\left(  x\right)  =\left\{
\begin{array}
[c]{cc}%
\varphi_{1}\left(  x-c\right)  & x\geq c\\
-\lambda\varphi_{1}\left(  -(x-c)\right)  & x<c
\end{array}
\right.  \label{U0}%
\end{equation}
where $\lambda\geq1$. Then (\ref{LAdefn}) is satisfied if and only if
$\lambda>1$, which renders $\varphi$ nondifferentiable at $c$. Thus a kink is
necessary for loss aversion given (\ref{U0}), but not given
(\ref{phi-trinomial}).

We add some interpretation of the functional form (\ref{phi-trinomial}). Take
$c=0$ for simplicity. Then, as observed above, loss aversion implies%
\[
\left(  x,\tfrac{1}{2};-x,\tfrac{1}{2}\right)  \prec0.
\]
How might one measure the degree of loss aversion expressed thereby? One
possibility is to use the reduction in the loss needed to imply indifference,
but then the new lottery would have nonzero mean which would obfuscate the
determination of "greater sensitivity to losses". Similarly if one were to
increase the odds of a gain with prizes unchanged. Thus we adjust both so as
to keep the zero mean. Specifically, we look for $\lambda>1$ such that
\begin{equation}
\left(  x,\lambda p;-\lambda x,p\right)  \sim0\text{ for all }x>0\text{ and
}0<p<1\text{.} \label{measure}%
\end{equation}
(For probabilities to sum to 1, one needs $p\left(  1+\lambda\right)  =1$, but
that can be safely ignored for present purposes given expected utility
theory.) The above condition states that when both the odds of a gain and the
size of the loss are increased by the factor $\lambda$, then (the zero mean
condition is satisfied and) indifference with $0$ is restored. In contrast,
when $\lambda=1$, then the strictly inferior $\tfrac{1}{2}/\tfrac{1}{2}$
lottery is obtained. This suggests using $\lambda-1$ to measure loss aversion.
Such a measure is well-defined for our model, using (\ref{phi-trinomial}),
since (\ref{measure}) is satisfied (uniquely) with $\lambda=\theta^{-1}$. Thus
$\theta^{-1}-1$ gives a measure of loss aversion that is \emph{behavioral}
(defined by the preference condition (\ref{measure})), and \emph{global} (the
same $\lambda$ works for all $x$ and $p$ as indicated).\footnote{In fact,
existence of $\lambda$ satisfying (\ref{measure}) is \emph{equivalent} to our
specification with $\lambda=\theta^{-1}$. More generally, one might weaken
(\ref{measure}) by allowing $\lambda$ to depend on $x$ and/or $p$. From that
perspective, our model yields a constant measure of loss aversion, perhaps
suggesting a partial analogue to CARA utility functions.} Alternatively, in
our model $\left(  \alpha x,p;-x,\alpha p\right)  \sim0$ is satisfied
(uniquely) by $\alpha=\theta<1$, suggesting $1-\theta$ as a measure of loss
aversion. In either case, the parameter $\theta$ admits a simple behavioral interpretation.

The results below, and the CLT underlying them, are limited to the case
$\theta=\underline{\sigma}/\overline{\sigma}$. However, they are robust to the
specification of $\varphi_{1}$, which is unrestricted except for nonparametric
monotonicity and concavity assumptions and technical (smoothness and
boundedness) conditions. In particular, what follows makes no assumption about
the nature of risk aversion in the domain of gains or alternatively, or about
the nature of risk loving in the domain of losses. The only relevant
restriction, imposed by (\ref{phi-trinomial}) and expected utility theory, is
on preference over \textquotedblleft mixed" lotteries.
%In addition to the interpretation of  given above, we can add the following behavioral implication
%dealing with mixed lotteries: If  is such that the indifference in (<ref>measure</ref>) is satisfied for some (x,p) pair,
%then it is satisfied for every other pair (x',p').

\subsubsection{Optimization\label{section-optimization}}

The preceding leads finally to the optimization problem (for each $n$)%
\begin{equation}
V_{n}\equiv\sup_{s\in\mathcal{S}}E_{P^{s}}\left[  \varphi\left(  \Sigma
_{1}^{n}X_{i}/\sqrt{n}\right)  \right]  \text{.} \label{Vn}%
\end{equation}
(Recall that $X_{i}\left(  \omega\right)  =\omega_{i}$ gives the outcome at
stage $i$.)

The finite horizon problem is not tractable (for us). For reasons of
tractability, Bayesian models in the literature typically take $\varphi$ to be
linear. In that case, given the fixed finite horizon $n$, the $\frac{1}%
{\sqrt{n}}$ factor is irrelevant and the objective is to maximize the expected
value of the sum $\Sigma_{1}^{n}X_{i}$. If outcomes are monetary prizes that
are perfect substitutes, which is the way we think of our model, then a linear
$\varphi$ implies risk neutrality as remarked in the introduction. An
alternative is that outcomes are measured in utils, as in the common
expected-additive-utility model of preference over risky consumption streams.
Then the underlying prizes (consumption levels, for example) at different
stages are not perfect substitutes, and also the ranking of the risky
consumption at stage $i$ is independent of the risks involved at other stages
(implying indifference to correlation in consumption risks). In applications
where these features are appropriate, indifference to risk in consumption (or
other underlying prizes) is not implied by a linear $\varphi$. However, for
the settings we have in mind, \emph{tractability comes at the cost of assuming
risk neutrality}.

Consider briefly a common approach to solving bandit problems analytically
which is to establish the optimality of index-based strategies, most commonly
using the Gittins index (Gittins and Jones 1974). When arms can be valued
separately, then at each stage and history an index summarizes each arm and
comparison of these indices determines which arm to pull. This approach does
not work in our model because arms cannot be delinked for at least two
reasons: (i) outcomes from one arm may be informative about the distribution
describing other arms because of common unknown parameters (see section
\ref{section-2arms}); (ii) because of loss aversion risk attitude depends on
the sign of the sum of past payoffs from \emph{all} arms.

Our approach to analysing (\ref{Vn}) for the loss averse utility index
(\ref{phi-trinomial}) is to study large-horizon approximations to the value
(indirect utility) of the bandit problem and corresponding approximately
optimal strategies. More precisely, define, conditional on showing below that
the following limit exists,
\begin{equation}
V\equiv\lim_{n\rightarrow\infty}V_{n}\text{.} \label{Vlim}%
\end{equation}
Below we derive results for $V$, which therefore imply approximate results for
$V_{n}$ when $n$ is sufficiently large. Secondly, say that the strategy
$s^{\ast}$ is \emph{asymptotically optimal} if
\begin{equation}
\lim_{n\rightarrow\infty}U_{n}\left(  s^{\ast}\right)  =\lim_{n\rightarrow
\infty}V_{n}\text{;} \label{ao}%
\end{equation}
or, equivalently, if, for every $\epsilon>0$, there exists $n^{\ast}$ such
that
\[
\mid U_{n}\left(  s^{\ast}\right)  -V_{n}\mid<\epsilon\text{ \ if }n>n^{\ast
}\text{.}%
\]
Thus asymptotic optimality of $s^{\ast}$ is a more concise way to say that
"$s^{\ast}$ is approximately optimal for problems with sufficiently long
horizon."\footnote{An implication is that, for any $s$, $\lim_{n\rightarrow
\infty}U_{n}\left(  s\right)  \leq\lim_{n\rightarrow\infty}U_{n}\left(
s^{\ast}\right)  $. This follows from (\ref{ao}) and $U_{n}\left(  s\right)
\leq V_{n}$ for all $n$.}

\subsection{\noindent Results\label{section-results}}

In all our results for the bandits model, the assumptions specified above are
adopted: conditional beliefs satisfy full support, (\ref{muzero1}) and
(\ref{condvar1}), and the utility index $\varphi$ is given by
(\ref{phi-trinomial}) and satisfies Assumption-Utility. Though the latter
requires $\overline{\sigma}>\underline{\sigma}$, all the results that follow
are trivially valid, by the classic martingale CLT, also when $\overline
{\sigma}=\underline{\sigma}$. Then all arms have a common variance and are
equivalent in the large horizon limit, making the (asymptotic) choice between
arms trivial. It simplifies discussions below to exclude that case.

\subsubsection{Value\label{section-value}}

Our first result concerns the limiting value $V$. We emphasize the surprising
(to us) degree to which this result is \emph{robust} to specifications of
$\varphi_{1}$ and the primitives $\left\{  P_{i}^{s_{i}}\right\}
_{i\geq1,~s_{i}\in\mathcal{S}_{i}}$, and therefore also to assumptions about
the nature of learning .

\begin{thm}
\label{thm-V}(i) Let $V_{n}$ be the value of the $n$-horizon problem
(\ref{Vn}). Then $\lim_{n\longrightarrow\infty}V_{n}$ exists. Moreover,
\begin{equation}
V=\lim_{n\longrightarrow\infty}V_{n}=\int_{-\infty}^{\infty}\varphi\left(
y\right)  q\left(  y\right)  dy\text{,} \label{Vq}%
\end{equation}
where $q$ is the pdf in (\ref{pdf-c>0})-(\ref{pdf-c<0}), which, for $c=0$
yields the simple form
\begin{equation}
q\left(  y\right)  =\left\{
\begin{array}
[c]{ccc}%
q^{\ast}\left(  y;\underline{\sigma}\right)  \left[  \frac{2\overline{\sigma}%
}{\underline{\sigma}+\overline{\sigma}}\right]  &  & y\geq0\\
&  & \\
q^{\ast}\left(  y;\overline{\sigma}\right)  \left[  \frac{2\underline{\sigma}%
}{\underline{\sigma}+\overline{\sigma}}\right]  &  & y<0
\end{array}
\right.  \label{qdensity}%
\end{equation}
Here $q^{\ast}(y;\sigma)=\frac{1}{\sqrt{2\pi}\sigma}\exp\left(  -({y/\sigma
})^{2}{/2}\right)  $ is the pdf for $\mathbb{N}\left(  0,\sigma^{2}\right)
$.\newline(ii) Let primitive beliefs be modified to $\left\{  \widehat
{P_{i}^{s_{i}}}\right\}  _{i\geq1,~s_{i}\in\mathcal{S}_{i}}$, another set
satisfying our assumptions, including counterparts of (\ref{muzero1}) and
(\ref{condvar1}), and where the latter is satisfied by the identical variance
extremes $\underline{\sigma}^{2}$ and $\overline{\sigma}^{2}$. Then
$\widehat{V}=V$.\newline(iii) The limiting value $V$ satisfies
\begin{equation}
V=\left\{
\begin{array}
[c]{cccc}%
= & \varphi(0) &  & c=0\\
> & \varphi(0) &  & c>0\\
< & \varphi(0) &  & c<0
\end{array}
\right.  \label{V=0}%
\end{equation}

\end{thm}

(i) not only proves that the large-horizon limit $V$ is well-defined, but also
gives an explicit description of $V$. Moreover, for some functions $\varphi$
the integral in (\ref{Vq}) can be expressed in closed form yielding a closed
form expression for $V$ for each $c$. For example, if $\varphi$ is taken to be
the exponential example (\ref{exp}), then, using the density in (\ref{pdf-c>0}%
)-(\ref{pdf-c<0}),%
\begin{equation}
V=\left\{
\begin{array}
[c]{ll}%
\!\!\Phi(-\frac{c}{\underline{\sigma}})-\Phi(\frac{c}{\underline{\sigma}%
})+e^{\frac{\underline{\sigma}^{2}}{2}}\left(  e^{-c}\Phi(-\underline{\sigma
}+\frac{c}{\underline{\sigma}})-e^{c}\Phi(-\underline{\sigma}-\frac
{c}{\underline{\sigma}})\right)  & \!\!c\leq0\\
\frac{\overline{\sigma}}{\underline{\sigma}}\left[  \Phi(-\frac{c}%
{\overline{\sigma}})-\Phi(\frac{c}{\overline{\sigma}})+e^{\frac{\underline
{\sigma}^{2}}{2}}\left(  e^{-\frac{\underline{\sigma}c}{\overline{\sigma}}%
}\Phi(-\underline{\sigma}+\frac{c}{\overline{\sigma}})-e^{\frac{\underline
{\sigma}c}{\overline{\sigma}}}\Phi(-\underline{\sigma}-\frac{c}{\overline
{\sigma}})\right)  \right]  & \!\!c>0,
\end{array}
\right.  \label{Vexp}%
\end{equation}
where $\Phi$ is the standard normal cdf.

The density $q$ in (\ref{qdensity}) yields a zero mean and variance equal to
$\underline{\sigma}\overline{\sigma}$, the geometric average of the two
extreme variances. Incorporation of the low (high) variance normal density for
positive (negative) arguments reflects risk aversion and loving on the two
subdomains respectively. Evidently, $q$ reduces to the normal density if
$\underline{\sigma}=\overline{\sigma}$, for example, there is a single arm.
Then (\ref{Vq}) is an immediate implication of the classic CLT. In the same
way, (i) follows directly from the new CLT in section \ref{section-CLT}.
Moreover, (i) is the main content of the theorem - the other parts follow
immediately from it. (ii) follows by inspection of the density and (iii)
follows from a simple calculation (see details in Appendix \ref{app-bandits}).

Part (ii) supports our hypothesis that the long-horizon heuristic reduces the
cognitive burden of the decision-maker. She need only know the variances of
arms, and even then, only for arms that have extreme variances.\footnote{The
identity of the extreme arms can vary with history. Therefore, there may not
exist two arms that are adequate for the entire horizon.} Here is some rough
intuition: Let the horizon be $n$ and consider the choice of arm at the last
stage given past realizations $x_{i}$ of $X_{i}$, $i<n$. It can be thought of
as maximizing $E_{P_{n}^{s}}\left[  \varphi\left(  ((\Sigma_{1}^{n-1}%
x_{i}+X_{n})/\sqrt{n}\right)  \right]  $ by choice of $s_{n}$ ($P_{n}^{s}$ is
the 1-step-ahead conditional in (\ref{agree})). The incremental payoff
$X_{n}/\sqrt{n}$ is small if $n$ is large. Thus a second-order Taylor series
expansion in $X_{n}$ can be used to approximate the objective function,
implying that the latter can be approximated (for each $s$) by a linear
function of both the mean (equal to zero by (\ref{muzero1})) and the
(conditional) variance. Finally, maximization of a linear function of variance
is necessarily achieved at an arm associated with either \underline{$\sigma$}
or $\overline{\sigma}$.
%Suppose make once-and-for-all choice of arm ex ante. Here also only variances matter for value, due to classic CLT.
%Moreover, under loss aversion the integral (EU) wrt normal dstribution is decreasing in variance. So optimal choice is
%necessarily to choose lowest variance.
%Aside: Does loss aversion imply that always reject a mean-preserving spread?

To interpret (iii), consider first the case $c=0$. Thus, for large $n$,
maximum expected utility is approximately equal to that achievable when the
payoff to each action is riskless, hence identically equal to the common mean,
implying zero gains and losses for sure. In other words, \emph{risk is a
matter of indifference in the limit}.
%The classical counterpart is where there is only one action, (equal sigmas), and phi(x)=--phi(-x).
%Then CLT implies that utility is 0 in limit because expected value of phi wrt Normal is zero, by symmetryof phi.
%Also, equality with 0 is "trivial" only if one takes our CLT as given!
%It is a corollary of the CLT and should be viewed as justifying its interest.
The freedom to switch between arms in response to experience is critical. If
one arm must be chosen ex ante for all trials, then maximum expected utility
is negative, hence less than $\varphi(0)=0$. (The classic CLT applies to each
arm separately and, by loss aversion, $\varphi\left(  -x\right)
<-\varphi\left(  x\right)  $ for all $x>0$; hence $\varphi\left(
\cdot\right)  $ has negative expected value under the normal $\mathbb{N}%
\left(  0,\sigma^{2}\right)  $ for any positive variance.) For further
perspective, consider the following lottery: Toss a fair coin. If Heads, then
receive a positive prize according to $\mathbb{N}\left(  0,\underline{\sigma
}^{2}\right)  $ conditioned on $\mathbb{R}_{+}$ and if Tails receive a
negative prize according to $\mathbb{N}\left(  0,\overline{\sigma}^{2}\right)
$ conditioned on $\mathbb{R}_{-}$. This lottery has negative expected utility
using $\varphi$. It is less attractive because the ability to choose actions
sequentially affords some influence over positive versus negative outcomes,
while in the lottery that influence belongs to nature alone.

Finally, (iii) implies that,\textit{ in the limit }$n\rightarrow\infty
$\textit{, a decision-maker with a positive reference point (}$c>0$)
\textit{strictly prefers the risky sequential choice problem to receiving zero
gain/loss for sure}.\textit{ }The intuition is that zero for sure is a certain
loss relative to a positive reference point, which makes it unattractive. A
positive reference point $c$ also reduces the limit value $V$, because it
reduces all gains and increases all losses ($\varphi\left(  x\right)
\overset{c}{\searrow}$ for all $x$), but to a lesser degree because of the
flexibility afforded by switching actions. Similarly, a negative reference
point implies the preference for the certain zero outcome. In this sense, a
higher benchmark or aspiration level leads to more participation in risky endeavors.

\begin{remark}
\label{remark-LLN}Suppose that DM uses the unweighted arithmetic average and
maximizes $E_{P^{s}}\left[  \varphi\left(  (\Sigma_{1}^{n}X_{i})/n\right)
\right]  $. Then a LLN would replace the CLT underlying (\ref{Vq}) and would
yield, by the LLN in Peng (2019, Theorem 2.4.1),
\begin{equation}
\lim_{n\rightarrow\infty}V_{n}=\varphi\left(  0\right)  =0\text{.}
\label{LLNlimit}%
\end{equation}
To reflect, consider the special case where there is independence across
trials of a single arm and across arms. Then by the classic LLN, the expected
utility of playing any $a\in\mathcal{A}$ at every stage and history converges
to $0\,$\ as $n\rightarrow\infty$. Consequently, for large $n$, DM is
approximately indifferent between repeated plays of $a$ and repeated plays of
any other $a^{\prime}$, because their means are identical. The implication of
(\ref{LLNlimit}) is that all such single-arm strategies are asymptotically
optimal, from which we conclude that, (in our setting, where only variances
differ), the LLN cannot serve as the basis for usefully approximating optimal
strategies for finite horizon problems. Furthermore, under the LLN,
(\ref{LLNlimit}) is valid not only for the loss averse functions $\varphi$
that we assume throughout, but also for all (suitably bounded and continuous)
$\varphi$ satisfying $\varphi\left(  0\right)  =0$. In contrast, in our model
using the $\sqrt{n}$-weighted average, such asymptotic risk neutrality is
satisfied only in the knife-edge case $c=0$, and risk is even strictly
desirable for $c>0$.
\end{remark}

\subsubsection{\medskip Strategies and the absence of
learning\label{section-nolearning}}

We describe an asymptotically optimal strategy for the special case where
there is no learning. The latter corresponds to the following restriction on
the primitive conditionals $\left\{  P_{i}^{s_{i}}\right\}  _{i\geq1,~s_{i}%
\in\mathcal{S}_{i}}$: For all $i\geq1$, $s_{i}\in\mathcal{S}_{i}$ and
histories $\left(  a^{\left(  i-1\right)  },\omega^{\left(  i-1\right)
}\right)  $,
\begin{equation}
P_{i}^{s_{i}}\left(  \cdot\mid a^{\left(  i-1\right)  },\omega^{(i-1)}\right)
=P_{1}^{s_{1}}\text{ \ if }s_{i}\left(  a^{\left(  i-1\right)  }%
,\omega^{\left(  i-1\right)  }\right)  =s_{1}\text{.} \label{nolearning}%
\end{equation}
Recall that at stage 1, history is null. Thus $s_{1}$ is simply an action and
$P_{1}^{s_{1}}$ gives (unconditional or) prior beliefs about the outcome of
action $s_{1}$. Thus (\ref{nolearning}) stipulates that for each given action
($s_{1}$ above), subsequent \emph{beliefs about the next outcome of that
action do not change with history} (where history includes past outcomes
associated with any, possibly different, action). An implication is that for
each fixed arm $a$, the joint probability distribution over outcomes given
repeated choice of $a$ is i.i.d. However, for other strategies $s$, the
induced measure $P^{s}$ (recall (\ref{Ps})) need not be a product measure.
(For example, if $\omega_{1}$ and $\omega_{1}^{\prime}$ are distinct outcomes,
and if $s$ specifies different actions at the histories $\left(  a_{1}%
,\omega_{1}\right)  $ and $(a_{1},\omega_{1}^{\prime})$, then the two
conditional probability distributions for stage $2$ outcomes generally differ.
This reflects a difference in the choice of action at stage 2 rather than
updating or learning.)

Define
\[
\sigma_{a}^{2}=E_{P_{1}^{s_{1}}}\left[  X_{1}^{2}\right]  \text{, \ if }%
s_{1}=a\in\mathcal{A}\text{.}%
\]
Then%
\[
\overline{\sigma}=\max_{a\in\mathcal{A}}\sigma_{a}\text{ and }\underline
{\sigma}=\min_{a\in\mathcal{A}}\sigma_{a}\text{.}%
\]

For simplicity, we focus first on $c=0$ and then indicate at the end of this
subsection how to accommodate $c\not =0$.

\begin{thm}
\label{thm-nolearn}Let $c=0$. Define strategy $s^{\ast}$ by $s_{1}^{\ast
}=\overline{a}$ and, for $n>1$,
\begin{equation}
s_{n}^{\ast}=\left\{
\begin{array}
[c]{ccc}%
\overline{a} & \text{if} & \Sigma_{1}^{n-1}X_{i}\leq0\\
\underline{a} & \text{if} & \Sigma_{1}^{n-1}X_{i}>0
\end{array}
\right.  \label{s*}%
\end{equation}
where $\sigma_{\overline{a}}=\overline{\sigma}$ and $\sigma_{\underline{a}}%
=$\underline{$\sigma$}. Then: (i) $s^{\ast}$ is asymptotically
optimal.\newline(ii) For every $N>0$,%
\begin{align*}
P^{s^{\ast}}\left(  \cap_{n=N}^{\infty}\{\Sigma_{1}^{n}X_{i}\leq0\}\right)
&  \leq\frac{\overline{\sigma}}{\overline{\sigma}+\underline{\sigma}}<1\text{
\ and }\\
P^{s^{\ast}}\left(  \cap_{n=N}^{\infty}\{\Sigma_{1}^{n}X_{i}>0\}\right)   &
\leq\frac{\overline{\sigma}}{\overline{\sigma}+\underline{\sigma}}<1.
\end{align*}
(iii) The high variance action is chosen less frequently in the limit. In
fact,
\begin{equation}
\lim_{n\rightarrow\infty}\frac{P^{s^{\ast}}\left(  \sigma_{s_{n}^{\ast}%
}=\overline{\sigma}\right)  }{P^{s^{\ast}}\left(  \sigma_{s_{n}^{\ast}%
}=\underline{\sigma}\right)  }=\frac{\underline{\sigma}}{\overline{\sigma}}<1.
\label{relativefreq}%
\end{equation}

\end{thm}

(i) identifies an asymptotically optimal $s^{\ast}$, while (ii) states that
$s^{\ast}$ exhibits switching between actions indefinitely with positive
probability according to the measure $P^{s^{\ast}}$ induced by $s^{\ast}$.
\ The latter fact indicates a difference between our model with loss aversion
and many bandit models. Commonly in the bandit literature, learning (or
exploration) provides the reason for switching, and eventually it is decided
that one arm is superior and experimentation ceases. Here, in contrast,
switching is optimal even in the absence of learning and (with positive
probability) persists indefinitely. This is because loss aversion implies that
the identity of the more attractive action or arm depends on whether one is in
a region of cumulative gains ($\Sigma_{1}^{n}X_{i}>0$) or cumulative losses
($\Sigma_{1}^{n}X_{i}<0$).\footnote{A global risk averter would choose the low
variance action $\underline{a}$ at every stage.} Finally, (iii) gives
explicitly the limiting relative frequencies induced by $s^{\ast}$.

We have emphasized the robustness of our results to the specification of
beliefs. An implication such robustness is given by comparing the present
no-learning model with a more general model where learning is admitted. Take
beliefs to be objective and consider DM being offered the choice between two
regimes, and suppose she chooses the one with the higher limiting value. One
regime corresponds to the conditional probability laws assumed to satisfy the
assumptions in section \ref{section-beliefs}. This regime accommodates many
different correlation patterns across trials. The other regime also satisfies
(\ref{nolearning}). Suppose further that the initial probability laws
$\{P_{1}^{s_{1}}\}_{s_{1}\in\mathcal{A}}$ are common to the two regimes. Then
these two specifications imply the same variance bounds $\overline{\sigma}$
and \underline{$\sigma$}. Therefore, by Theorem \ref{thm-V}, they imply the
same limiting value $V$.\ Consequently, \emph{the correlation (or lack
thereof) between outcomes in distinct trials of an action is a matter of
indifference in the limit}.\footnote{In the case of a single arm, the above
reduces to the following familiar property of the classic martingale CLT. Let
$P\in\Delta\left(  \Pi_{1}^{\infty}\Omega_{i}\right)  $ satisfy the
counterparts of (\ref{muzero1}) and (\ref{condvar1}). Then they are also
satisfied by $Q$, the i.i.d. product of the marginal $P_{1}\in\Delta\left(
\Omega_{1}\right)  $. Moreover, $P$ and $Q$ imply the same limiting
probability distribution for $\Sigma_{1}^{n}X_{i}/\sqrt{n}$.}

\begin{remark}
It does not follow that there is an asymptotically optimal strategy common to
both regimes. Theorem \ref{thm-learn} in the next section describes an optimal
strategy under learning for a special case.
\end{remark}

Finally, we describe how the theorem can be extended to accommodate $c\not =%
0$. For that purpose, instead of using a single strategy to approximate
finite-horizon problems, consider a \emph{sequence }$s^{n}=\left(  s_{i}%
^{n}\right)  $\emph{ of strategies}, where, for each $n$, $s^{n}\in
\mathcal{S}$ is thought of as a strategy used in the $n$-horizon problem
(\ref{Vn}). (Accordingly, components $s_{i}^{n}$ with $i>n$ are irrelevant.)
The counterpart of (\ref{ao}) is
\begin{equation}
\lim_{n\rightarrow\infty}E_{P^{s^{n}}}[\varphi\left(  \Sigma_{1}^{n}%
X_{i}/\sqrt{n}\right)  ]=\lim_{n\rightarrow\infty}V_{n}=V \label{ao-n}%
\end{equation}
Then, arguing as in the proof of Theorem \ref{thm-nolearn}, one can show that
(\ref{ao-n}) is satisfied by $s^{n}$, where, for each $n\geq1$ and $1\leq
i\leq n$,{\textbf{ }}%
\[
s_{i}^{n}=\left\{
\begin{array}
[c]{ccc}%
\overline{a} & \text{if} & \Sigma_{1}^{i-1}X_{j}/\sqrt{n}\leq c\\
\underline{a} & \text{if} & \Sigma_{1}^{i-1}X_{j}/\sqrt{n}>c\text{.}%
\end{array}
\right.
\]
$s_{i}^{n}$ can be defined arbitrarily if either $n=1$ or $i>n$.

\subsubsection{A classic two-armed bandit problem
revisited\label{section-2arms}}

There are two arms, $a$ and $b$, hence $\mathcal{A}=\{a,b\}$. Thus the set of
possible outcomes for each arm and stage is $\overline{\Omega}=\{1,-1,0\}$,
and outcomes are governed, both ex ante and \emph{for any history}, by the
following probabilities:
\begin{align*}
\text{ arm }a  &  \text{: }\Pr\left(  1\right)  =\Pr\left(  -1\right)
=p_{a}/2\\
\text{arm }b  &  \text{: }\Pr\left(  1\right)  =\Pr\left(  -1\right)
=p_{b}/2.
\end{align*}
For each arm, outcomes follow a random walk with zero mean and with variance
equal to the appropriate value of $p$. It is known that
\begin{equation}
\{p_{a},p_{b}\}=\{\underline{p},\overline{p}\}\text{,} \label{pbars}%
\end{equation}
where $0<\underline{p}<\overline{p}<1$ are known; equivalently, the set of
variances $\{\underline{\sigma}^{2},\overline{\sigma}^{2}\}$ is known, where
\[
\underline{\sigma}^{2}=\underline{p}\text{ and }\overline{\sigma}%
^{2}=\overline{p}.
\]
However, there is uncertainty about which of $\underline{p}$ and $\overline
{p}$ describes arm $a$ and which describes arm $b$, that is, there is
uncertainty about which arm has the higher variance. DM has prior beliefs
about which arm is which, and forms Bayesian posteriors as experience
accumulates. At each stage, she chooses which arm to pull taking into account
what she has learned about the arms from past experience.

\begin{remark}
\label{remark-2armed}Uncertainty about "which arm is which" in a 2-arm setting
is a classic version of the bandit problem (Bradt, Johnson and Karlin 1956;
Feldman 1962); indeed, the former refer to it (p. 1060) as \textquotedblleft%
\underline{the} Two-armed Bandit." These and subsequent papers typically
assume a finite horizon and maximization of the expected value of the sum of
payoffs, (in particular, means rather than variances are the focus).
\end{remark}

Our framework accommodates the above learning process. The set of primitive
conditionals $\{P_{n}^{s_{n}}\}_{n\geq1,s_{n}\in\mathcal{S}_{n}}$ is defined
as follows. DM's prior beliefs about which arm is which are completely
specified by $\mu_{1}$, the probability she assigns initially to
$p_{a}=\underline{p}$. Thus, prior probabilities of the outcomes from choosing
arm $\alpha$, $\alpha=a,b$, are given by%
\begin{align*}
P_{1}^{a}\left(  1\right)   &  =\mu_{1}\underline{p}/2+(1-\mu_{1})\overline
{p}/2=P_{1}^{a}\left(  -1\right) \\
P_{1}^{b}\left(  1\right)   &  =(1-\mu_{1})\underline{p}/2+\mu_{1}\overline
{p}/2=P_{1}^{b}\left(  -1\right)  \text{,}%
\end{align*}
which can be expressed in terms of our formalism by{
\begin{align*}
P_{1}^{s_{1}}\left(  \omega_{1}\right)   &  =I_{\{s_{1}=a,\omega_{1}\neq
0\}}\left[  \mu_{1}\underline{p}/2+(1-\mu_{1})\overline{p}/2\right] \\
&  +I_{\{s_{1}=a,\omega_{1}=0\}}\left[  \mu_{1}(1-\underline{p})+(1-\mu
_{1})(1-\overline{p})\right] \\
&  +I_{\{s_{1}=b,\omega_{1}\neq0\}}\left[  (1-\mu_{1})\underline{p}/2+\mu
_{1}\overline{p}/2\right] \\
&  +I_{\{s_{1}=b,\omega_{1}=0\}}\left[  (1-\mu_{1})(1-\underline{p})+\mu
_{1}(1-\overline{p})\right]  \text{.}%
\end{align*}
For later stages, }DM updates her prior probability that $p_{a}=\underline{p}$
to the Bayesian posterior $\mu_{n}$, $n>1$, defined inductively by%
\begin{align}
&  \log\left(  \frac{\mu_{n+1}/(1-\mu_{n+1})}{\mu_{n}/(1-\mu_{n})}\right)
\label{posteriors}\\
&  =\left[  I_{a}\left(  a_{n}\right)  -I_{b}\left(  a_{n}\right)  \right]
\left(  \left(  1-I_{0}\left(  \omega_{n}\right)  \right)  \log\left(
\frac{\underline{p}}{\overline{p}}\right)  +I_{0}\left(  \omega_{n}\right)
\log\left(  \frac{1-\underline{p}}{1-\overline{p}}\right)  \right)
\text{.}\nonumber
\end{align}
Then the conditional probability $P_{n}^{s_{n}}$, for each $n>1$ and stage
strategy $s_{n}$, { is given by%
\begin{align}
P_{n}^{s_{n}}\left(  \omega_{n}|a^{(n-1)},\omega^{(n-1)}\right)   &
=I_{\{s_{n}=a,\omega_{n}\neq0\}}\left[  \mu_{n}\underline{p}/2+(1-\mu
_{n})\overline{p}/2\right] \nonumber\\
&  +I_{\{s_{n}=a,\omega_{n}=0\}}\left[  \mu_{n}(1-\underline{p})+(1-\mu
_{n})(1-\overline{p})\right] \label{Psn}\\
&  +I_{\{s_{n}=b,\omega_{n}\neq0\}}\left[  (1-\mu_{n})\underline{p}/2+\mu
_{n}\overline{p}/2\right] \nonumber\\
&  +I_{\{s_{n}=b,\omega_{n}=0\}}\left[  (1-\mu_{n})(1-\underline{p})+\mu
_{n}(1-\overline{p})\right]  \text{.}\nonumber
\end{align}
}

Consider also the probability measure $P^{s}$, for $s\in\mathcal{S}$,
constructed as in (\ref{Ps}) by pasting the above conditionals. It is
completely described by its restriction to finite dimensional cylinders, and
thus view $P^{s}$ as a measure on $\Pi_{1}^{n}\Omega_{i}$. For any
$\omega^{\left(  n\right)  }=\left(  \omega_{1},...,\omega_{n}\right)  $, the
outcomes of the first $n$ trials, and the given $s$, define the induced
frequency vector $f^{s}\left(  \omega^{\left(  n\right)  }\right)  $,
\begin{equation}
f^{s}\left(  \omega^{\left(  n\right)  }\right)  =(f_{a}^{s}\left(
\omega^{\left(  n\right)  }),f_{b}^{s}\left(  \omega^{\left(  n\right)
}\right)  ,f_{a,0}^{s}(\omega^{\left(  n\right)  }),f_{b,0}^{s}\left(
\omega^{\left(  n\right)  }\right)  \right)  \text{,} \label{f}%
\end{equation}
where: for $\alpha\in\{a,b\}$, $\ f_{\alpha}^{s}\left(  \omega^{\left(
n\right)  }\right)  $ and $f_{\alpha,0}^{s}\left(  \omega^{\left(  n\right)
}\right)  $ give, respectively, the number of trials of arm $\alpha$ and the
number of those that yield the outcome $0$. Then the ex ante probability of
the above outcomes are given by\footnote{The proof is elementary and is
omitted.}%
\begin{align}
P^{s}\left(  \omega_{1},...,\omega_{n}\right)   &  =\mu_{1}\left[
(\underline{p}/2)^{f_{a}^{s}-f_{a,0}^{s}}\left(  \overline{p}/2\right)
^{f_{b}^{s}-f_{b,0}^{s}}\left(  1-\underline{p}\right)  ^{f_{a,0}^{s}}\left(
1-\overline{p}\right)  ^{f_{b,0}^{s}}\right] \label{Probn}\\
&  +(1-\mu_{1})\left[  \left(  \overline{p}/2\right)  ^{f_{a}^{s}-f_{a,0}^{s}%
}\left(  \underline{p}/2\right)  ^{f_{b}^{s}-f_{b,0}^{s}}\left(
1-\overline{p}\right)  ^{f_{a,0}^{s}}\left(  1-\underline{p}\right)
^{f_{b,0}^{s}}\right]  \text{.}\nonumber
\end{align}
The two terms on the right correspond to the two possible scenarios, $p_{a}%
=$\underline{$p$} or $\overline{p}$, weighted by their prior probabilities.
Conditional on each scenario the expression reflects two assumptions: (i)
independence between distinct trials, whether conducted with the same arm or
with different arms; and (ii) all trials with a given arm are viewed as
similar (or interchangeable) so that the probability of any (finite) sequence
of outcomes for that arm is invariant to any reordering (accordingly, for each
arm, the probability of a set of outcomes depends only on the number of
occurrences of $0$ and $\{1,-1\}$). This latter assumption of "symmetry"
within each arm is known as \emph{partial exchangeability,} a property\emph{
}introduced by\emph{ }de Finetti (1938), who also showed that it implies
conditional independence as in (i), and, in fact, that it characterizes a
representation such as in (\ref{Probn}).\footnote{The stronger property of
exchangeability, which is better known, assumes interchangeability also across
distinct arms and thus views the two arms as being identical, which is
excluded in our case because of (\ref{pbars}) and $\underline{p}%
\not =\overline{p}$. See Link (1980) and Diaconis and Freedman (1982) for more
on partial exchangeability and Kallenberg (2005) for a comprehensive treatment
of probabilistic symmetries.}

The preceding satisfies all the assumptions of our general model and hence
Theorem \ref{thm-V} applies. Moreover, with the added structure assumed
herein\ we can also address strategies and what is learned asymptotically.
Below we assume $c=0$.

Define the strategy $s^{\ast}$ by $s_{1}^{\ast}=a$ and, for $n>1$,%
\[
s_{n}^{\ast}=\left\{
\begin{array}
[c]{ccc}%
a & \text{if} &
\begin{array}
[c]{c}%
\Sigma_{1}^{n-1}X_{j}\leq0\text{, }\mu_{n}<\frac{1}{2}\text{ OR}\\
\!\Sigma_{1}^{n-1}X_{j}>0\text{, }\mu_{n}>\frac{1}{2}%
\end{array}
\\
&  & \\
b & \text{if} & \text{otherwise}%
\end{array}
\right.
\]
According to $s^{\ast}$, arm $a$ is used at stage $n>1$ if (and only if) there
are cumulative losses and it is more likely that $a$ has higher variance
($\mu_{n}<\frac{1}{2}$), or there are cumulative gains and it is more likely
that $a$ has lower variance ($\mu_{n}>\frac{1}{2}$). Intuition argues for this
choice of arm at stage $n$ if there are no later trials remaining, but may
seem myopic more generally. Nevertheless, we show that $s^{\ast}$ is
approximately optimal for large horizons. (For other instances where myopic
strategies are optimal in bandit problems see, for example, Banks and Sundaram
(1992) and the papers cited in Remark \ref{remark-2armed}.)

\begin{thm}
\label{thm-learn} Let $c=0$ and $\mu_{1}\in\lbrack0,1]$. Then $s^{\ast}$ is
asymptotically optimal.
\end{thm}

\noindent When $\mu_{1}\in\{0,1\}$, we are back in the no-learning case of the
last section and Theorem \ref{thm-nolearn} applies.

Conclude with observations about the process of posteriors $\{\mu_{n}\}$ that
confirm for our setting properties familiar from Bayesian learning theory.

\begin{remark}
\label{remark-posteriors}Let $s\in\mathcal{S}$ be any strategy. Then:\newline%
(i) Posteriors converge to certainty, that is, for any prior $\mu_{1}%
$,\footnote{In fact, convergence to certainty is valid for every $P^{s}$,
$s\in\mathcal{S}.$}
\begin{equation}
\lim_{n\rightarrow\infty}\mu_{n}\in\{0,1\}\text{ \ }P^{s}\text{-a.s. }
\label{mulimit}%
\end{equation}
(ii) Suppose that, unknown to the decision-maker, the truth is that $p_{a}%
=$\underline{$p$}. Consequently, given any strategy $s$, outcomes are governed
by the probability law $Q^{s}\in\Delta\left(  \Pi_{1}^{\infty}\Omega
_{i},\mathcal{G}\right)  $, whose $1$-step-ahead conditionals are $Q_{i}^{s}$,
$i\geq1$, given by%
\[
Q_{i}^{s}\left(  1\right)  =Q_{i}^{s}\left(  -1\right)  =\left\{
\begin{array}
[c]{ccc}%
\underline{p}/2 & \text{if} & s_{i}=a\\
\overline{p}/2 & \text{if} & s_{i}=b
\end{array}
\right.
\]
Then\textbf{,} for every $\mu_{1}>0$,
\begin{equation}
\lim_{n\rightarrow\infty}\mu_{n}=1\text{ \ }Q^{s}\text{-a.s.} \label{mulimit2}%
\end{equation}

\end{remark}

Think of $\{\mu_{n}\}$ as representing subjective beliefs. Then (\ref{mulimit}%
) expresses the decision-maker's ex ante complete confidence that
asymptotically she will know "which arm is which." In (ii), $Q^{s}$ is the
true probability law over outcome sequences when strategy $s$ is adopted, and
hence (\ref{mulimit2}) is an expression of "Bayesian consistency". Both
results are valid for any strategy, and thus reflect Bayesian updating alone
and not asymptotic optimality.

\section{A Central Limit Theorem}

\subsection{Preliminaries}

The mathematical basis for our analysis of the bandit problem is a central
limit theorem about sets of measures that will be provided here. To smooth the
transition for the reader, we begin with a few remarks about connect the
bandit model to sets of measures.

In section \ref{section-beliefs}, we introduced the primitive set of
one-step-ahead conditionals $\left\{  P_{i}^{s_{i}}\right\}  _{i\geq
1,~s_{i}\in\mathcal{S}_{i}}$, and then pointed out that, for each $s=\left(
s_{1},...,s_{i},...\right)  $, these conditionals can be pasted together to
obtain a measure $P^{s}\in\Delta\left(  \Pi_{1}^{\infty}\Omega_{i}%
,\mathcal{G}\right)  $. Now we collect all these measures and define the set
$\mathcal{P}\subset\Delta\left(  \Pi_{1}^{\infty}\Omega_{i},\mathcal{G}%
\right)  $ by
\begin{equation}
\mathcal{P}=\{P^{s}:s\in\mathcal{S}\}. \label{Pbandit}%
\end{equation}
Our CLT will be applied to this set. However, in order to better reveal its
underlying structure and to facilitate other potential applications, (for
example, to models concerned with robustness to model uncertainty), the CLT
will be formulated and proven more generally. In particular, rather than
restricting ourselves to the sets associated with the bandit problem, the CLT
will take as a primitive an abstract set $\mathcal{P}\subset\Delta\left(
\Pi_{1}^{\infty}\Omega_{i},\mathcal{G}\right)  $ satisfying "rectangularity,"
defined in the next section, which we suggest is key to the CLT, and which we
show (Lemma \ref{lemma-banditrectang}) is satisfied by the set defined in
(\ref{Pbandit}).

One more observation is helpful for the transition. For the set $\mathcal{P}$
defined by (\ref{Pbandit}), it is immediate that, for each $n$,
\begin{equation}
V_{n}=\sup_{s\in\mathcal{S}}E_{P^{s}}[\varphi\left(  \Sigma_{1}^{n}X_{i}%
/\sqrt{n}\right)  ]=\sup\limits_{Q\in\mathcal{P}}E_{Q}[\varphi\left(
\Sigma_{1}^{n}X_{i}/\sqrt{n}\right)  ]\text{.} \label{equiv}%
\end{equation}
{\normalsize \hfill}

\noindent The CLT will involve expressions such as that on the right in
(\ref{equiv}). However, we can think of the supremum over measures as
equivalent to optimization over strategies.

\medskip

To proceed, adopt the mathematical primitives $\left(  \Pi_{1}^{\infty}%
\Omega_{i},\{\mathcal{G}_{n}\}_{n=1}^{\infty}\right)  $ and $\mathcal{G}$,
though with possibly different interpretations.\footnote{In fact, we do not
need the previous assumptions that $\Omega_{i}$ is identical for all $i$ and
finite. Here the $\Omega_{i}$s are arbitrary.} For each $i\geq1$, $X_{i}:$
$\Pi_{1}^{\infty}\Omega_{j}\longrightarrow\mathbb{R}$ is $\mathcal{G}_{i}%
$-measurable. Another primitive is a set $\mathcal{P}\subset\Delta\left(
\Pi_{1}^{\infty}\Omega_{i},\mathcal{G}\right)  $, not to be confused with the
set in (\ref{Pbandit}). The following assumptions are adopted throughout: All
measures in $\mathcal{P}$ are equivalent on each $\mathcal{G}_{n}$,
\begin{equation}
E_{Q}[X_{i}|\mathcal{G}_{i-1}]=0\text{ for all }Q\in\mathcal{P}\text{ and all
}i\geq1\text{,} \label{muzero}%
\end{equation}
and conditional variances satisfy, for some $\overline{\sigma}\geq
\underline{\sigma}>0$,
\begin{equation}
ess\sup\limits_{Q\in\mathcal{P}}E_{Q}\left[  X_{i}^{2}|\mathcal{G}%
_{i-1}\right]  =\overline{\sigma}^{2}\text{ and }ess\inf_{Q\in\mathcal{P}%
}E_{Q}\left[  X_{i}^{2}|\mathcal{G}_{i-1}\right]  =\underline{\sigma}%
^{2}\text{ for all }i\geq1\text{.} \label{condvar}%
\end{equation}
Assume also that $\left(  X_{i}\right)  $ satisfies the \emph{Lindeberg
condition}:
\begin{equation}
\lim_{n\rightarrow\infty}\frac{1}{n}\sum\limits_{i=1}^{n}\sup\limits_{Q\in
\mathcal{P}}E_{Q}\left[  \left\vert X_{i}\right\vert ^{2}I_{\{\left\vert
X_{i}\right\vert >\sqrt{n}\epsilon\}}\right]  =0,\quad\forall\epsilon>0.
\label{linder}%
\end{equation}

\noindent When $\mathcal{P}$ is a singleton, these conditions reduce to
standard assumptions in the classic martingale CLT.

The following notation is convenient in the sequel: $\mathcal{H}$ denotes the
set of all random variables $X$ on $\left(  \Pi_{1}^{\infty}\Omega
_{i},\mathcal{G}\right)  $ satisfying $\sup_{Q\in\mathcal{P}}E_{Q}[\left\vert
X\right\vert ]<\infty$. For any $X$ in $\mathcal{H}$, its (conditional)
\emph{upper expectations} are defined by
\[
\mathbb{E}[X]\equiv\sup_{Q\in\mathcal{P}}E_{Q}[X]\text{,}\;\ \mathbb{E}\left[
X|\mathcal{G}_{n}\right]  \equiv ess\sup\limits_{Q\in\mathcal{P}}E_{Q}\left[
X|\mathcal{G}_{n}\right]  \text{.}%
\]
\hfill\hfill

\subsection{Rectangularity\label{section-rectangularity}}

The final assumption underlying our CLT is that the set $\mathcal{P}$ is
"rectangular", which means that it is closed with respect to the pasting of
alien 1-step-ahead conditionals. Rectangularity was introduced in Epstein and
Schneider (2003) in the context of recursive utility theory, where an
axiomatic analysis demonstrated its role in modeling dynamic behavior for an
ambiguity-averse decision-maker.\footnote{It has been studied and applied also
in robust stochastic dynamic optimization (Iyengar, 2005; Shapiro, 2016), in
the literature on dynamic risk measures (Riedel, 2004; Cheridito, Delbaen and
Kupper, 2006; Acciaio and Penner, 2011), and in continuous-time modeling in
finance (Chen and Epstein, 2002).} When $\mathcal{P}$ is the singleton
$\{P\}$, rectangularity is trivially implied by Bayesian updating,
specifically by the fact that after decomposing $P$ into its 1-step-ahead
conditionals, these can be pasted together to recover $P$. More generally,
rectangularity requires that the set $\mathcal{P}$ is closed also with respect
to pasting together 1-step-ahead conditionals that are alien, that is, induced
by possibly different measures in $\mathcal{P}$. This property endows
$\mathcal{P}$ with a recursive structure that yields a form of the law of
iterated expectations (see below).

For a formal definition, we introduce some additional notation. Write
\begin{align*}
\omega_{\left(  n\right)  }  &  =\left(  \omega_{n},...\right)  \text{,
}\omega^{\left(  n\right)  }=\left(  \omega_{1},...,\omega_{n}\right)  ,\\
\mathcal{P}_{0,n}  &  =\{P_{\mid\mathcal{G}_{n}}:P\in\mathcal{P}\}\text{
and}\\
\mathcal{G}_{(n+1)}  &  =\left\{  A\subset\Pi_{n+1}^{\infty}\Omega_{i}:\Pi
_{1}^{n}\Omega_{i}\times A\in\mathcal{G}\right\}  \text{.}%
\end{align*}
A probability kernel from $\left(  \Pi_{1}^{n}\Omega_{i},\mathcal{G}%
_{n}\right)  $ to $\left(  \Pi_{n+1}^{\infty}\Omega_{i},\mathcal{G}%
_{(n+1)}\right)  $ is a function $\lambda:\Pi_{1}^{n}\Omega_{i}\times
\mathcal{G}_{(n+1)}\longrightarrow\left[  0,1\right]  $ satisfying:

\begin{description}
\item[Kernel 1:] $\forall\omega^{(n)}\in\Pi_{1}^{n}\Omega_{i}$, $\lambda
\left(  \omega^{\left(  n\right)  },\cdot\right)  $ is a probability measure
on $\left(  \Pi_{n+1}^{\infty}\Omega_{i},\mathcal{G}_{(n+1)}\right)  $,

\item[Kernel 2:] $\forall A\in\mathcal{G}_{(n+1)}$, $\lambda\left(
\cdot,A\right)  $ is a $\mathcal{G}_{n}$-measurable function on $\Pi_{1}%
^{n}\Omega_{i}$.
\end{description}

Any pair $\left(  p_{n},\lambda\right)  $ consisting of a probability measure
$p_{n}$ on $\left(  \Pi_{1}^{n}\Omega_{i},\mathcal{G}_{n}\right)  $ and a
probability kernel $\lambda$ as above, induces a unique probability measure
$P$ on $\left(  \Pi_{1}^{\infty}\Omega_{i},\mathcal{G}\right)  $ that
coincides with $p_{n}$ on $\mathcal{G}_{n}$. It is given by, $\forall
A\in\mathcal{G}$,%
\begin{equation}
P\left(  A\right)  =\int\limits_{\Pi_{1}^{n}\Omega_{i}}\int\limits_{\Pi
_{n+1}^{\infty}\Omega_{i}}I_{A}\left(  \omega^{\left(  n\right)  }%
,\omega_{(n+1)}\right)  \lambda\left(  \omega^{\left(  n\right)  }%
,d\omega_{(n+1)}\right)  p_{n}\left(  d\omega^{\left(  n\right)  }\right)
\text{.} \label{Pn}%
\end{equation}

For $Q\in\mathcal{P}$, let $Q\left(  \cdot\mid\mathcal{G}_{n}\right)  $,
denote its induced (regular) conditional. Then it defines a probability kernel
$\lambda$ by: $\forall\omega^{\left(  n\right)  }\in\Pi_{1}^{n}\Omega_{i}$,
\begin{equation}
\lambda\left(  \omega^{\left(  n\right)  },A\right)  =Q\left(  \Pi_{1}%
^{n}\Omega_{i}\times A\mid\mathcal{G}_{n}\right)  \left(  \omega^{\left(
n\right)  }\right)  \text{, }\forall A\in\mathcal{G}_{(n+1)}\text{. ~}
\label{kernel}%
\end{equation}

\noindent A feature of such a kernel is that the single measure $Q$ is used to
define the conditional at every $\omega^{\left(  n\right)  }$. We are
interested in kernels for which the measure to be conditioned can vary with
$\omega^{\left(  n\right)  }$. Thus say that the probability kernel $\lambda$
is a $\mathcal{P}$-\emph{kernel} if: $\forall\omega^{\left(  n\right)  }\in
\Pi_{1}^{n}\Omega_{i}$ $\exists Q\in\mathcal{P}$ such that (\ref{kernel}) is satisfied.

Finally, say that $\mathcal{P}$ is \emph{rectangular }(with respect to the
filtration $\{\mathcal{G}_{n}\}$) if: $\forall n~\forall p_{n}\in$
$\mathcal{P}_{0,n}$ and for every $\mathcal{P}$-\emph{kernel }$\lambda$, if
$P$ is defined as in (\ref{Pn}), then $P\in\mathcal{P}$.

The significance of rectangularity is illuminated by the following lemma. (Its
proof can be found in Chen and Epstein (2020).)

\begin{lemma}
{\normalsize \label{lemma-rect} $\mathcal{P}$ }rectangular implies the
following (for any $0\leq m\leq n\in N$).

\noindent(i) \textbf{Stability by composition:} For any $Q,R\in\mathcal{P}$,
$\exists P\in\mathcal{P}$ such that, for any $X\in\mathcal{H}$,
\[
E_{P}[X|\mathcal{G}_{m}]=E_{Q}[E_{R}[X|\mathcal{G}_{n}]|\mathcal{G}%
_{m}]\text{.}%
\]
(ii) \textbf{Stability by bifurcation:} For any $Q,R\in\mathcal{P}$, and any
$A_{n}\in\mathcal{G}_{n}$, $\exists P\in\mathcal{P}$ such that, for any
$X\in\mathcal{H}$,
\[
E_{P}[X|\mathcal{G}_{n}]=I_{A_{n}}E_{Q}[X|\mathcal{G}_{n}]+I_{A_{n}^{c}}%
E_{R}[X|\mathcal{G}_{n}].
\]
\noindent(iii) \textbf{Law of iterated upper expectations}: For any
$X\in\mathcal{H}$,
\begin{equation}
\mathbb{E}[\mathbb{E}[X|\mathcal{G}_{n}]|\mathcal{G}_{m}]=\mathbb{E}%
[X|\mathcal{G}_{m}]. \label{LIE}%
\end{equation}
(iv) Let $\{X_{i}\}$ be a sequence in $\mathcal{H}$.
%Set $S_{n-1}=\sum_{i=1}^{n-1}X_{i}$.
Then, for any continuous bounded functions $f,h$:{
\[
\mathbb{E}\left[  f\left(  \Sigma_{1}^{n-1}X_{i}\right)  +h\left(  \Sigma
_{1}^{n-1}X_{i}\right)  X_{n}^{2}\right]  =\mathbb{E}\left[  \mathbb{E}\left[
f\left(  \Sigma_{1}^{n-1}X_{i}\right)  +h\left(  \Sigma_{1}^{n-1}X_{i}\right)
X_{n}^{2}|\mathcal{G}_{n-1}\right]  \right]  .
\]
} If the conditional means and variances of $X_{n}$ satisfy (\ref{muzero}) and
(\ref{condvar}), then
\[
\mathbb{E}\left[  h\left(  \Sigma_{1}^{n-1}X_{i}\right)  X_{n}^{2}%
|\mathcal{G}_{n-1}\right]  =\overline{\sigma}^{2}\left[  h(\Sigma_{1}%
^{n-1}X_{i})\right]  ^{+}-\underline{\sigma}^{2}\left[  h(\Sigma_{1}%
^{n-1}X_{i})\right]  ^{-}.
\]
(Superscripts $+$ and $-$ denote the positive and negative parts respectively.)
\end{lemma}

\noindent

Part (iii) gives the law of iterated expectations for upper expectations (a
similar condition for lower expectations is implied). (iv) is an extension
that is used in the proofs of our CLTs. Parts (i) and (ii) of the lemma
describe direct implications of $\mathcal{P}$ being "closed with respect to
the pasting of alien marginals and conditionals."

Importantly, the set of measures generated in the bandit model is rectangular.

\begin{lemma}
\label{lemma-banditrectang}The set $\mathcal{P}$ defined in (\ref{Pbandit}) is rectangular.
\end{lemma}

\noindent\textbf{Proof: }{Because of the finiteness of each $\Omega_{i}$ in
the bandit model, it suffices to consider the following. Let $P^{s}$,
$P^{s^{\prime}}$ and $P^{s^{\prime\prime}}$ be measures in $\mathcal{P}$. Take
$n>1$ and $A_{n}\in\mathcal{G}_{n}$, and paste the measures together by
constructing the new measure $P$ by,
\[
P\left(  B\right)  =\int\limits_{\Pi_{1}^{n}\Omega_{i}}\int\limits_{\Pi
_{n+1}^{\infty}\Omega_{i}}I_{B}\left(  \omega^{\left(  n\right)  }%
,\omega_{(n+1)}\right)  \lambda\left(  \omega^{\left(  n\right)  }%
,d\omega_{(n+1)}\right)  P_{|\mathcal{G}_{n}}^{s}\left(  d\omega^{\left(
n\right)  }\right)  \text{,}\ \forall B\in\mathcal{G},
\]
where, for all $(\omega^{\left(  n\right)  },A)\in\Pi_{i=1}^{n}\Omega
_{i}\times\mathcal{G}_{(n+1)}$,%
\[
\lambda\left(  \omega^{\left(  n\right)  },A\right)  =I_{A_{n}}\left(
\omega^{\left(  n\right)  }\right)  P^{s^{\prime}}\left(  \Pi_{1}^{n}%
\Omega_{i}\times A\mid\mathcal{G}_{n}\right)  +I_{A_{n}^{c}}\left(
\omega^{\left(  n\right)  }\right)  P^{s^{\prime\prime}}\left(  \Pi_{1}%
^{n}\Omega_{i}\times A\mid\mathcal{G}_{n}\right)  \text{.}%
\]
}Then $P=P^{\widehat{s}}$, where $\widehat{s}=\left(  s_{1},...,s_{n}%
,\widehat{s}_{n+1},\widehat{s}_{n+2},...\right)  \in\mathcal{S}$ is defined
by, for any $i>n$ and $\omega^{(i-1)}=\left(  \omega^{(n)},\cdots,\omega
_{i-1}\right)  ,$
\begin{align*}
\widehat{s}_{i}\left(  \omega^{\left(  i-1\right)  }\right)   &  =I_{A_{n}%
}\left(  \omega^{\left(  n\right)  }\right)  s_{i}^{\prime}\left(  a^{\prime
}{}^{(i-1)},\omega^{\left(  i-1\right)  }\right)  +I_{A_{n}^{c}}\left(
\omega^{\left(  n\right)  }\right)  s_{i}^{\prime\prime}\left(  a^{\prime
\prime}{}^{(i-1)},\omega^{\left(  i-1\right)  }\right)  \text{, }\\
\text{and }a_{j}^{\prime}  &  =s_{j}^{\prime}\left(  a^{\prime}{}%
^{(j-1)},\omega^{\left(  j-1\right)  }\right)  \text{,}\ a_{j}^{^{\prime
\prime}}=s_{j}^{\prime\prime}\left(  a^{\prime\prime}{}^{(j-1)},\omega
^{\left(  j-1\right)  }\right)  \text{, for }1\leq j\leq
i-1.\ \ \ \blacksquare
\end{align*}
\bigskip

As a concrete example, consider again the special case of the bandit model
with no-learning (section \ref{section-nolearning}). In that model, for each
arm $a$, the probability distribution over outcomes on a single trial is
$P_{1}^{a}$, ($P_{1}^{a}\equiv$ $P_{1}^{s_{1}}$ where $s_{1}=a$), independent
of history. Denote by $\mathcal{L}$ the set of all such measures over outcomes
as $a$ varies over all arms. Since any arm can be chosen at any history, the
corresponding set of measures over outcome sequences is\footnote{Formally, it
follows from (\ref{nolearning}) and (\ref{agree}).}
\[
\mathcal{P}=\left\{  P\in\Delta\left(  \Pi_{1}^{\infty}\Omega_{i}%
,\mathcal{G}\right)  :P_{i}\in\mathcal{L}\text{ for every }i\text{ and
history}\right\}  \text{.}%
\]
($P_{i}$ is the 1-step-ahead conditional at stage $i$ induced by $P$.) This
set is obviously rectangular. Note that in spite of $\mathcal{L}$ being common
to all trials and histories, trial outcomes as modeled by $\mathcal{P}$ are
not necessarily identical. Indeed, any measure in $\mathcal{L}$ can describe
the $i^{th}$ trial at a specific history in conjunction with any possibly
different measure in $\mathcal{L}$ being the law describing the $j^{th}$ trial
at any other history. As a result, besides the restriction imposed by
$\mathcal{L}$, the set $\mathcal{P}$ imposes no restrictions on the pattern of
heterogeneity across trials.
%Just restating that any action is feasible at any history.  But perhaps worthwhile to make clear that
%rectangulairty for a set set of measures  captures unstructured heterogeneity.
In particular, it \emph{accommodates heteroscedasticity that is largely
unstructured} apart from the restrictions imposed by the extreme variances
\underline{$\sigma$} and $\overline{\sigma}$. To varying degrees, the same is
true for all rectangular sets associated with our bandit model.

\subsection{The theorem\label{section-CLT}}

We extend (a version of) the classic martingale CLT to admit a set of
variances while maintaining the assumption of a fixed zero mean. Throughout
$(B_{t})$ denotes a standard Brownian motion under a probability space
$(\Omega^{\ast},\mathcal{F}^{\ast},P^{\ast})$ and $(\mathcal{F}_{t})_{t\geq0}$
is the natural filtration generated by $(B_{t})$.

In the classic case, the limiting distribution is normal, which is the
distribution of $B_{1}$. In the more general case, the corresponding (upper)
limit is not given by the normal distribution, but is described instead by the
time 1 value of an \emph{oscillating Brownian motion} (Keilson and Wellner
1978; Lejay and Pigato 2018), defined as follows: Given $\overline{\sigma}%
\geq\underline{\sigma}>0$ and threshold $c\in\mathbb{R}$, let $(W_{t}^{c})$
denote the unique strong solution, (which exists by Le Gall (1984)), of the
stochastic differential equation (SDE)%
\begin{equation}
Y_{t}=\int_{0}^{t}\sigma\left(  Y_{s}\right)  dB_{s},\quad t\geq0,
\label{sde0}%
\end{equation}
where the diffusion coefficient $\sigma$ is the positive two-valued function,
discontinuous at the threshold $c$,
\begin{equation}
\sigma(y)=\underline{\sigma}I_{[c,\infty)}(y)+\overline{\sigma}I_{(-\infty
,c)}(y),\ \forall y\in\mathbb{R}. \label{sigma-function}%
\end{equation}
There is a seeming connection to the bandit model - lower volatility in the
region $\left(  c,\infty\right)  $ of gains where there is risk aversion, and
large volatility in the region of losses $\left(  -\infty,c\right)  $ where
there is risk loving.\footnote{Reversing the roles of $\overline{\sigma}$ and
\underline{$\sigma$} also defines an oscillating Brownian motion, but one that
is irrelevant here given the assumption of loss aversion.} In fact, by Keilson
and Wellner (1978, Theorem 1), the time 1 value $W_{1}^{c}$ of the oscillating
Brownian motion has distribution given by the density $q$ referred to in
Theorem \ref{thm-V}(i).

\begin{thm}
\label{thm-CLT}
%Let the sequence $\left(  X_{i}\right)  $ be such
%that $X_{i}\in\mathcal{H}$ for each $i$, and where $\left(  X_{i}\right)  $
%satisfies (\ref{mubar}) and (\ref{condvar}), with zero  conditional
%means, and conditional upper and lower variances $\overline{\sigma}^2\ge\underline{\sigma}^{2}>0$. Assume also the Lindeberg condition (\ref{linder})
%and that $\mathcal{P}$ is rectangular. Set $S_n:=\sum_{i=1}^nX_{i}$, then,
Let the sequence $\left(  X_{i}\right)  $ be such that $X_{i}\in\mathcal{H}$
for each $i$, and where $\left(  X_{i}\right)  $ satisfies (\ref{muzero}) and
(\ref{condvar}), with conditional upper and lower variances $\overline{\sigma
}^{2}\geq\underline{\sigma}^{2}>0$. Assume also the Lindeberg condition
(\ref{linder}), that measures in $\mathcal{P}$ are equivalent on each
$\mathcal{G}_{i}$, and that $\mathcal{P}$ is rectangular. Set $\theta
=\underline{\sigma}/\overline{\sigma}$. For any $c\in\mathbb{R}$ and
$\varphi_{1}\in C_{b}^{3}\left(  \mathbb{R}_{+}\right)  $, with $\varphi
_{1}\left(  0\right)  =0$, define $\varphi$ by
\begin{equation}
\varphi(x)=\left\{
\begin{array}
[c]{lll}%
\ \ \varphi_{1}(x-c) &  & x\geq c\\
-\frac{1}{\theta}\varphi_{1}\left(  -\theta(x-c)\right)  &  & x<c
\end{array}
\right.  \label{function-phbar}%
\end{equation}

\noindent If $\varphi_{1}^{\prime\prime}(x)\leq0$ for $x\geq0$,
%we define
%$$\varphi_2(x+c)=-\frac{\overline{\sigma}}{\underline{\sigma}}\varphi_1\left(-\frac{\underline{\sigma}}{\overline{\sigma}}x+c\right)+\left(1+\frac{\overline{\sigma}}{\underline{\sigma}}\right)\varphi_1(c),\ \forall x\le0,$$
%and
%\begin{equation*}
%\varphi(x)=\left\{
%\begin{array}{ll}
%\varphi_{1}(x),& x\ge c;\\
%\varphi_{2}(x)=-\frac{\overline{\sigma}}{\underline{\sigma}}\varphi_1\left(-\frac{\underline{\sigma}}{\overline{\sigma}}(x-c)+c\right)+\left(1+\frac{\overline{\sigma}}{\underline{\sigma}}\right)\varphi_1(c),& x<c.
%\end{array}
%\right.
%\end{equation*}
then
\begin{equation}
\lim_{n\rightarrow\infty}\sup\limits_{Q\in\mathcal{P}}E_{Q}\left[
\varphi\left(  \frac{\Sigma_{1}^{n}X_{i}}{\sqrt{n}}\right)  \right]
=E_{P^{\ast}}[\varphi(W_{1}^{c})]. \label{thm-clt-eq-O-leq-sup}%
\end{equation}

\end{thm}

%\begin{remark}
%A question is can we find any functions $\varphi$ that has no dependence on $(\underline{\sigma},\overline{\sigma})$ such the Theorem \ref{CLT-O} holds.
%Let $\tilde{\varphi}\in C_b^3(\mathbb{R})$ be a positive homogeneity odd function, that is
%$$\tilde{\varphi}(x)=-\tilde{\varphi}(-x)\text{ and }\tilde{\varphi}(\lambda x)=\lambda\tilde{\varphi}(x),\ \forall \lambda\ge0.$$
%If $\tilde{\varphi}''(x)\ge0$ for $x\ge0$ (or  $\tilde{\varphi}''(x)\le0$ for $x\ge0$), Theorem \ref{CLT-O} still holds for  $\tilde{\varphi}$.
%Also, for the indicator function for the interval   $(-\infty,c]$ or $[c,\infty)$ the above Theorem still holds, more details in Corollary \ref{cor-pro-O}.
%\end{remark}
\
%If c=0, then the RHS in each of the 4 equations equals phi(0).
%Thus oscillating BM does not play a role, except in proof, if c=0. It is relevant, however, when we admit nonzero c.
%Further, even given c=0, the equality to phi(0) is not obvious intuitively, (symmetry arises only in the limit?),
%and an alternative proof is not clear.
%In bandits section, c=0 is the main hypothesis, but we do make some use also of CLT when c is nonzero.

The most important point to make about the theorem is that all its assumptions
are satisfied by the bandit model with $\mathcal{P}$ defined by (\ref{Pbandit}%
). (The Lindeberg condition (\ref{linder}) is satisfied because of the
finiteness of $\Omega_{i}=\overline{\Omega}$.) Therefore, using also the noted
density for $W_{1}^{c}$, the CLT implies Theorem \ref{thm-V}(i). Though the
bandit theorem is stated with reference only to a density and not to
oscillating Brownian motions, we prefer to include the latter here because it
is more revealing of what underlies the limit and, to a degree, how the limit
result is proven.

For perspective, if instead of defining $\varphi$ by (\ref{function-phbar}),
we took $\varphi$ to be any (suitably bounded, smooth and) \emph{globally
concave} function, then the limit in (\ref{thm-clt-eq-O-leq-sup}) would equal
the expected value of $\varphi$ under $\mathbb{N}\left(  0,\underline{\sigma
}^{2}\right)  $, as in the classic case with fixed variance $\underline
{\sigma}^{2}$. Informally, this result is suggested by taking $c\rightarrow
-\infty$ above. (For a rigorous argument, see Proposition 2.2.15 and Theorem
2.4.4 in Peng (2019).)

Some extensions of the CLT are possible. For example, one can obtain similar
limits with any combination of the modifications $\theta=\frac{\overline
{\sigma}}{\underline{\sigma}}$, $\varphi_{1}^{\prime\prime}(x)\geq0$ on
$\left(  0,\infty\right)  $, and/or one considers the limit of the lower
expectation $\inf_{Q\in\mathcal{P}}E_{Q}\left[  \varphi\left(  \Sigma_{1}%
^{n}X_{i}/\sqrt{n}\right)  \right]  $. These extensions do not seem relevant
to the bandit problem, but the reader can find them in our working paper
version listed in the bibliography.\textbf{ }It is also possible to derive
closed-form limiting results for other integrands (functions $\varphi$), for
example, for some indicator functions (Appendix \ref{app-indicators}). For
many other functions $\varphi$, the corresponding expressions for the limit
are more complex, less transparent and arguably intractable, and consequently
are excluded.

We conclude with mention of related CLTs in the literature. Chen and Epstein
(2020) establish CLTs assuming, contrary to (\ref{muzero})-(\ref{condvar}),
that conditional means lie in an interval $\left[  \underline{\mu}%
,\overline{\mu}\right]  $ while all conditional variances equal a constant
$\sigma^{2}$. In common with this paper, rectangularity is a key assumption.
However, their theorems are substantially different, for example, limits have
a different form and proofs are much different. There exist other
generalizations of the classic CLT that are motivated by robustness to
ambiguity. In both Marinacci (1999, Theorem 16) and Epstein, Kaido and Seo
(2016), experiments are not ordered and their analyses are better suited for a
cross-sectional, rather than sequential, context. Another difference is that
in both cases, limiting distributions are normal. Peng (2007, 2019) and Fang
et al (2019) assume that experiments are ordered. Comparison with Theorem 3.2
of the latter is representative. It is more general than our results, for
example, in permitting ambiguity about both mean and variance. For purposes of
comparison, limit attention to the special case of their theorem where there
is ambiguity about variance only. Even then, an important difference,
particularly given the application developed here, is that greater generality
comes arguably at the cost of reduced tractability. In particular, limits are
much more complicated (they involve Peng's (2007) notion of a "G-normal"
distribution), and a counterpart of Theorem \ref{thm-CLT} is not apparent from
their results.\footnote{Another difference, apart from applicability, is that
Peng and coauthors take a nonlinear expectation operator as the core primitive
and adopt the PDE approach, while our primitive is a set of probability
measures. Notably, conditionals are central in our analysis, including in the
notion of rectangularity, but are not easily accommodated in the operator
approach. Thus their key assumption is similar in spirit to rectangularity,
but is expressed in terms of nonlinear expectations and without recourse to
conditionals.} Finally, none of the above papers recognize the potential
application to sequential decision problems such as the bandit problem.

\appendix

\section{Appendix: Main Proofs\label{app-CLT}}

The notation and assumptions in Theorem \ref{thm-CLT} are adopted throughout
this appendix. Let $(B_{t})$ be the standard Brownian motion under a
probability space $(\Omega^{\ast},\mathcal{F}^{\ast},P^{\ast})$, and let
$(\mathcal{F}_{t})_{t\geq0}$ be the natural filtration generated by
$(B_{t})_{t\geq0}$.
%We often denote $\Sigma_{1}^{n}X_{i}$ by $S_{n}$.

\subsection{Lemmas}

{\normalsize For a small fixed $h>0$, and any fixed $(t,x,c)\in\lbrack
0,1+h]\times\mathbb{R}\times\mathbb{R}$, $(Y_{s}^{t,x,c})_{s\in\lbrack
t,1+h]}$ denotes the solution of the SDE
\begin{equation}
\left\{
\begin{array}
[c]{l}%
dY_{s}^{t,x,c}=\sigma\left(  Y_{s}^{t,x,c}\right)  dB_{s},\quad s\in\lbrack
t,1+h]\\
Y_{t}^{t,x,c}=x,
\end{array}
\right.  \label{sde3}%
\end{equation}
where $\sigma(y)=$\underline{$\sigma$}$I_{[c,\infty)}(y)+\overline{\sigma
}I_{(-\infty,c)}(y),\ \forall y\in\mathbb{R}$. }

{\normalsize By Keilson and Wellner (1978, Theorem 1), (see also Chen and Zili
(2015)), the transition probability density of $(Y_{s}^{t,x,c})_{s\in\lbrack
t,1+h]}$ is given by, for any $t<s\leq1+h$ and $y\in\mathbb{R}$},{\normalsize
\begin{align}
q^{c}(t,x;s,y)=  &  \frac{1}{\sqrt{2\pi(s-t)}}\frac{1}{\sigma(y)}\exp\left(
-\frac{\left(  \frac{x-c}{\sigma(x)}-\frac{y-c}{\sigma(y)}\right)  ^{2}%
}{2(s-t)}\right) \nonumber\\
&  +\frac{\overline{\sigma}-\underline{\sigma}}{\overline{\sigma}%
+\underline{\sigma}}\frac{1}{\sqrt{2\pi(s-t)}}\frac{sgn(y-c)}{\sigma(y)}%
\exp\left(  -\frac{\left(  \left\vert \frac{x-c}{\sigma(x)}\right\vert
+\left\vert \frac{y-c}{\sigma(y)}\right\vert \right)  ^{2}}{2(s-t)}\right)
.\quad\label{tranproba}%
\end{align}
}

{\normalsize
%We consider a special class of functions in $C_b^3(\mathbb{R})$, which is defined by
%$$\Psi(\alpha,\beta,c):=\{\varphi\in C_b^3(\mathbb{R})\ |\ \alpha\varphi''(\alpha y+c)=-\beta\varphi''(-\beta y+c),\ \forall y\geq0 \}.$$
Given $\varphi_{1}\in C_{b}^{3}(\mathbb{R}_{+})$, $\varphi$ is defined by
(\ref{function-phbar}). Then
\[
\varphi\in C_{b}^{1}(\mathbb{R})\text{ and }\varphi^{\prime\prime
}(z+c)=-\tfrac{\underline{\sigma}}{\overline{\sigma}}\varphi^{\prime\prime
}(-\tfrac{\underline{\sigma}}{\overline{\sigma}}z+c),\ \forall z<0.
\]
}

{\normalsize
}

{\normalsize
%%%%%%%%%%%%%%%%%%%%%%%%%%%%%%%%%%%%%%%%%%%%%%%%%%%%%%%%%%%%%%%%%%%%%%%%%%%%%%%%%%%%%%%%%%%%%%%%%%%%%%%%%%%%%%%%%%%%%%%%%%%%%%%%%%%%%%%%%
}

{\normalsize Define the set of functions $\{H_{t}\}_{t\in\lbrack0,1+h]}$ by
%Take time interval [0,1+h] only to justify 2nd order differentiability [in what?] at 1%
\begin{equation}
H_{t}(x)=E_{P^{\ast}}\left[  \varphi\left(  Y_{1+h}^{t,x,c}\right)  \right]
,\quad\forall x\in\mathbb{R}. \label{hf}%
\end{equation}
Then
\[
H_{1+h}(x)=\varphi(x),\quad H_{0}(0)=E_{P^{\ast}}[\varphi(Y_{1+h}%
^{0,0,c})]=E_{P^{\ast}}[\varphi(W_{1+h}^{c})].
\]
}

{\normalsize The following lemma describes some properties of the functions
$\{H_{t}\}_{t\in\lbrack0,1+h]}$. }

\begin{lemma}
{\normalsize \label{lemma-ddp}The functions $\{H_{t}\}$ defined by (\ref{hf})
satisfy: }

\begin{description}
\item[(1)] {\normalsize For any $t\in\lbrack0,1]$, $H_{t}\in C_{b}%
^{2}(\mathbb{R})$, and the first and second derivatives of $H_{t}$ are bounded
uniformly in $t\in\lbrack0,1]$. }

\item[(2)] {\normalsize There exists a constant $L$ such that, for any
$x_{1},x_{2}\in\mathbb{R}$ and $t\in\lbrack0,1]$,
\[
\left\vert H_{t}^{\prime\prime}(x_{1})-H_{t}^{\prime\prime}(x_{2})\right\vert
\leq L|x_{1}-x_{2}|.
\]
}

\item[(3)] If $\varphi^{\prime\prime}(x)\leq0$ for $x>c$,
%$${\bf( C2).}\quad\varphi''(x)\leq0\text{ for }x\geq c\text{ and }\varphi''(x)\geq0\text{ for }x\leq c,$$
then{\normalsize
\[
\left\{
\begin{array}
[c]{ll}%
H_{t}^{\prime\prime}(x)\leq0 & \text{for }x\geq c\\
H_{t}^{\prime\prime}(x)\geq0 & \text{for }x\leq c\text{.}%
\end{array}
\right.
\]
}

{\normalsize
%\color{red}
%{\bf Remark: the sign of $H_{t}^{\prime\prime}$ only depend on the sign of $\varphi''$ and independent with the value of $\alpha,\beta$.}
%\color{black}
}

\item[(4)] {\normalsize For any $r\in\lbrack0,1+h-t]$,
\[
H_{t}(x)=E_{P^{\ast}}\left[  H_{t+r}\left(  Y_{t+r}^{t,x,c}\right)  \right]
,\ \forall x\in\mathbb{R}.
\]
}

\item[(5)] {\normalsize If $\varphi^{\prime\prime}(x)\leq0$ for $x>c$, then
\[
\lim_{n\rightarrow\infty}\sum_{m=1}^{n}\sup\limits_{x\in\mathbb{R}}\left\vert
H_{\frac{m-1}{n}}\left(  x\right)  -H_{\frac{m}{n}}\left(  x\right)
-\frac{\overline{\sigma}^{2}}{2n}\left[  {H_{\frac{m}{n}}^{\prime\prime}%
(x)}\right]  ^{+}+\frac{\underline{\sigma}^{2}}{2n}\left[  {H_{\frac{m}{n}%
}^{\prime\prime}(x)}\right]  ^{-}\right\vert =0\text{.}%
\]
}

\item[(6)] {\normalsize There exists a constant $C_{0}$ such that
\[
\sup_{x\in\mathbb{R}}\left\vert H_{1}(x)-\varphi(x)\right\vert \leq C_{0}%
\sqrt{\underline{\sigma}^{2}+\overline{\sigma}^{2}}\sqrt{h}.
\]
}
\end{description}
\end{lemma}

{\normalsize \noindent\textbf{Proof:} \textbf{(1)} Given the transition
probability density in (\ref{tranproba}), we have, for }$t\in\left[
0,1\right]  $,{\normalsize
\[
H_{t}(x)=\int_{-\infty}^{\infty}\varphi(y)q^{c}(t,x;1+h,y)dy,\quad\forall
x\in\mathbb{R}.
\]
For $T=1+h$, we have
\begin{align*}
H_{t}^{\prime}(x)  &  =\left\{
\begin{array}
[c]{ll}%
\frac{1}{\underline{\sigma}\sqrt{2\pi(T-t)}} \int^{\infty}_{0} \varphi
_{1}^{\prime}(y) \left[  e^{-\frac{(x-c+y)^{2}}{2\underline{\sigma}^{2}(T-t)}%
}+e^{-\frac{(x-c-y)^{2}}{2\underline{\sigma}^{2}(T-t)}}\right]  dy & \text{if
}x\geq c\\
\frac{1}{\overline{\sigma}\sqrt{2\pi(T-t)}} \int_{0}^{\infty}\varphi
_{1}^{\prime}(\tfrac{\underline{\sigma}}{\overline{\sigma}}y) \left[
e^{-\frac{(x-c+y)^{2}}{2\overline{\sigma}^{2}(T-t)}}+e^{-\frac{(x-c-y)^{2}%
}{2\overline{\sigma}^{2}(T-t)}}\right]  dy & \text{ if }x\leq c
\end{array}
\right. \\
H_{t}^{\prime\prime}(x)  &  =\left\{
\begin{array}
[c]{ll}%
\frac{1}{\underline{\sigma}\sqrt{2\pi(T-t)}}\int_{0}^{\infty}\varphi
_{1}^{\prime\prime}(y)e^{-\frac{(x-c-y)^{2}}{2\underline{\sigma}^{2}(T-t)}%
}\left[  1-e^{-\frac{2y(x-c)}{\underline{\sigma}^{2}(T-t)}}\right]  dy &
\text{if }x\geq c\\
\frac{-1}{\overline{\sigma}\sqrt{2\pi(T-t)}}\int_{0}^{\infty}\tfrac
{\underline{\sigma}}{\overline{\sigma}}\varphi_{1}^{\prime\prime}%
(\tfrac{\underline{\sigma}}{\overline{\sigma}}y)e^{-\frac{(x-c+y)^{2}%
}{2\overline{\sigma}^{2}(T-t)}}\left[  1-e^{\frac{2y(x-c)}{\overline{\sigma
}^{2}(T-t)}}\right]  dy & \text{if }x\leq c\text{.}%
\end{array}
\right.
\end{align*}
%\color{red}
%Actually,
%$$
%H_{t}^{\prime\prime}(x)=\left\{
%\begin{array}{ll}
%\!\!\frac{1}{\alpha\sqrt{2\pi T}}\left[\frac{2\beta}{\beta+\alpha}(\varphi_1'(c)-\varphi_2'(c))e^{-\frac{(x-c)^2}{2\alpha^2T}}+\int_{0}^{\infty}\varphi_1''(y+c)e^{-\frac{(x-c-y)^2}{2\alpha^2T}}dy\right.&\\
%\hspace{1cm}\left.+\int_{-\infty}^0\left(\frac{\beta-\alpha}{\alpha+\beta}\varphi_1''(-y+c)+\frac{2\beta^2}{\alpha(\alpha+\beta)}\varphi_2''(\tfrac{\beta}{\alpha}y+c)\right)e^{-\frac{(x-c-y)^2}{2\alpha^2T}}dy\right],& x\geq c\\
%\!\!\frac{1}{\beta\sqrt{2\pi T}}\left[\frac{2\alpha}{\beta+\alpha}(\varphi_1'(c)-\varphi_2'(c))e^{-\frac{(x-c)^2}{2\beta^2T}}+\int_{-\infty}^0\varphi_2''(y+c)e^{-\frac{(x-c-y)^2}{2\beta^2T}}dy\right.&\\
%\hspace{1cm}\left.+\int_0^{\infty}\left(\frac{\alpha-\beta}{\alpha+\beta}\varphi_2''(-y+c)+\frac{2\alpha^2}{\beta(\alpha+\beta)}\varphi_1''(\tfrac{\alpha}{\beta}y+c)\right)e^{-\frac{(x-c-y)^2}{2\beta^2T}}dy\right],&x\leq c
%\end{array}
%\right.
%$$
%Since $\varphi_1'(c)=\varphi_2'(c)$ and $\varphi_2''(y+c)=-\frac{\alpha}{\beta}\varphi_1''(-\frac{\alpha}{\beta}y+c),\ y<c$, we can deduce (\ref{sec-de-ht}).
%\color{black}
}

\noindent The assertion follows from {\normalsize $\varphi_{1}\in C_{b}%
^{3}(\mathbb{R}_{+})$ and the definition of $\varphi$ in (\ref{function-phbar}%
). }

{\normalsize
%When $t=1$, $H_1(x)\equiv\varphi(x)$, it's trivial.
}

\smallskip

\noindent{\normalsize \textbf{(2)} For any $x>c$, \ }$H_{t}^{\prime
\prime\prime}(x)=${\small
\[
\frac{1}{\underline{\sigma}\sqrt{2\pi(T-t)}}\left[  2\varphi_{1}^{\prime
\prime}(0)e^{-\frac{(x-c)^{2}}{2\underline{\sigma}^{2}(T-t)}}+\int_{0}%
^{\infty}\varphi_{1}^{\prime\prime\prime}(y)\left(  e^{-\frac{(x-c-y)^{2}%
}{2\underline{\sigma}^{2}(T-t)}}+e^{-\frac{(x-c+y)^{2}}{2\underline{\sigma
}^{2}(T-t)}}\right)  dy\right]  ,
\]
}and{\small , }{\normalsize for $x<c$, \ }$H_{t}^{\prime\prime\prime}%
(x)=${\small
\[
\frac{1}{\overline{\sigma}\sqrt{2\pi(T-t)}}\left[  -2\varphi_{1}^{\prime
\prime}(0)e^{-\frac{(x-c)^{2}}{2\overline{\sigma}^{2}(T-t)}}+\int^{\infty}%
_{0}\tfrac{\underline{\sigma}^{2}}{\overline{\sigma}^{2}}\varphi_{1}%
^{\prime\prime\prime}(\tfrac{\underline{\sigma}}{\overline{\sigma}}y)\left(
e^{-\frac{(x-c-y)^{2}}{2\overline{\sigma}^{2}(T-t)}}+e^{-\frac{(x-c+y)^{2}%
}{2\overline{\sigma}^{2}(T-t)}}\right)  dy\right]  .
\]
}{\normalsize Since $\varphi_{1}\in C_{b}^{3}(\mathbb{R}_{+})$, there exists a
constant $L$ such that
\[
\sup\limits_{x\in\mathbb{R},x\neq c}|H_{t}^{\prime\prime\prime}(x)|\leq
L\text{ ~for all }t\in\left[  0,1\right]  \text{.}%
\]
The assertion follows by the Mean Value Theorem. }

\smallskip

\noindent{\normalsize \textbf{(3)} It follows from the explicit form of
$H_{t}^{\prime\prime}(x)$ given above. }

\smallskip

\noindent{\normalsize \textbf{(4)} Since $(Y_{s}^{t,x,c})$ is a
time-homogeneous Markov process, for any $r\in\lbrack0,1+h-t]$,
\[
H_{t}(x)=E_{P^{\ast}}[\varphi(Y_{1+h}^{t,x,c})]=E_{P^{\ast}}\left[
E_{P^{\ast}}[\varphi(Y_{1+h}^{t,x,c})|\mathcal{F}_{t+r}]\right]  =E_{P^{\ast}%
}[H_{t+r}(Y_{t+r}^{t,x,c})].
\]
}

\smallskip

\noindent{\normalsize \textbf{(5)} It follows from part (4) that, for any
$1\leq m\leq n$,%
\[
H_{\frac{m-1}{n}}(x)=E_{P^{\ast}}\left[  H_{\frac{m}{n}}\left(  Y_{\frac{m}%
{n}}^{\frac{m-1}{n},x,c}\right)  \right]  .
\]
Apply It\^{o}'s formula to $H_{\frac{m}{n}}\left(  Y_{\frac{m}{n}}^{\frac
{m-1}{n},x,c}\right)  $ to derive
\begin{align*}
H_{\frac{m}{n}}\left(  Y_{\frac{m}{n}}^{\frac{m-1}{n},x,c}\right)  =  &
H_{\frac{m}{n}}\left(  x\right)  +\int_{\frac{m-1}{n}}^{\frac{m}{n}}%
H_{\frac{m}{n}}^{\prime}\left(  Y_{s}^{\frac{m-1}{n},x,c}\right)
\sigma\left(  Y_{s}^{\frac{m-1}{n},x,c}\right)  dB_{s}\\
&  +\frac{1}{2}\int_{\frac{m-1}{n}}^{\frac{m}{n}}H_{\frac{m}{n}}^{\prime
\prime}\left(  Y_{s}^{\frac{m-1}{n},x,c}\right)  \left(  \sigma\left(
Y_{s}^{\frac{m-1}{n},x,c}\right)  \right)  ^{2}ds
\end{align*}
Using parts (3) and (4), we have
\begin{align*}
H_{\frac{m-1}{n}}\left(  x\right)   &  =E_{P^{\ast}}\left[  H_{\frac{m}{n}%
}\left(  Y_{\frac{m}{n}}^{\frac{m-1}{n},x,c}\right)  \right]  =\\
&  E_{P^{\ast}}\left[  H_{\frac{m}{n}}\left(  x\right)  +\frac{\overline
{\sigma}^{2}}{2}\int_{\frac{m-1}{n}}^{\frac{m}{n}}\left[  {H_{\frac{m}{n}%
}^{\prime\prime}\left(  Y_{s}^{\frac{m-1}{n},x,c}\right)  }\right]
^{+}ds-\frac{\underline{\sigma}^{2}}{2}\int_{\frac{m-1}{n}}^{\frac{m}{n}%
}\left[  {H_{\frac{m}{n}}^{\prime\prime}\left(  Y_{s}^{\frac{m-1}{n}%
,x,c}\right)  }\right]  ^{-}ds\right]
\end{align*}
Thus%
\begin{align*}
&  \sum_{m=1}^{n}\sup\limits_{x\in\mathbb{R}}\left\vert H_{\frac{m-1}{n}%
}\left(  x\right)  -H_{\frac{m}{n}}\left(  x\right)  -\frac{\overline{\sigma
}^{2}}{2n}\left[  {H_{\frac{m}{n}}^{\prime\prime}(x)}\right]  ^{+}%
+\frac{\underline{\sigma}^{2}}{2n}\left[  {H_{\frac{m}{n}}^{\prime\prime}%
(x)}\right]  ^{-}\right\vert \\
\leq &  \sum_{m=1}^{n}\sup\limits_{x\in\mathbb{R}}E_{P^{\ast}}\left[
\frac{\underline{\sigma}^{2}+\overline{\sigma}^{2}}{2}\int_{\frac{m-1}{n}%
}^{\frac{m}{n}}\left\vert H_{\frac{m}{n}}^{\prime\prime}\left(  Y_{s}%
^{\frac{m-1}{n},x,c}\right)  -H_{\frac{m}{n}}^{\prime\prime}(x)\right\vert
ds\right] \\
\leq &  \sum_{m=1}^{n}\sup\limits_{x\in\mathbb{R}}\frac{(\underline{\sigma
}^{2}+\overline{\sigma}^{2})L}{2n}E_{P^{\ast}}\left[  \sup_{s\in\lbrack
\frac{m-1}{n},\frac{m}{n}]}\left\vert Y_{s}^{\frac{m-1}{n},x,c}-x\right\vert
\right] \\
\leq &  \sum_{m=1}^{n}\sup\limits_{x\in\mathbb{R}}\frac{C}{n}\left(
E_{P^{\ast}}\left[  \int_{\frac{m-1}{n}}^{\frac{m}{n}}\left(  \sigma\left(
Y_{s}^{\frac{m-1}{n},x,c}\right)  \right)  ^{2}dr\right]  \right)  ^{\frac
{1}{2}}\leq\frac{C\sqrt{\underline{\sigma}^{2}+\overline{\sigma}^{2}}}%
{\sqrt{n}}\text{,}%
\end{align*}
where $C$ is a constant that depends only on \underline{$\sigma$}%
$,\overline{\sigma},L$. }

{\normalsize \smallskip}

{\normalsize \noindent}\noindent

\noindent{\normalsize \textbf{(6)} Since $\varphi\in C_{b}^{1}(\mathbb{R})$,
$C_{0}\equiv\Vert\varphi^{\prime}\Vert=\sup_{x\in\mathbb{R}}|\varphi^{\prime
}(x)|<\infty$, and%
\begin{align*}
\sup_{x\in\mathbb{R}}\left\vert H_{1}(x)-\varphi(x)\right\vert =  &
\sup_{x\in\mathbb{R}}\left\vert E_{P^{\ast}}[\varphi(Y_{1+h}^{1,x,c}%
)]-\varphi(x)\right\vert \\
\leq &  \sup_{x\in\mathbb{R}}E_{P^{\ast}}\left[  \left\vert \varphi
(Y_{1+h}^{1,x,c})-\varphi(x)\right\vert \right] \\
\leq &  \sup_{x\in\mathbb{R}}C_{0}E_{P^{\ast}}\left[  \left\vert \int
_{1}^{1+h}\sigma\left(  Y_{s}^{1,x,c}\right)  dB_{s}\right\vert \right] \\
\leq &  \sup_{x\in\mathbb{R}}C_{0}\left(  E_{P^{\ast}}\left[  \int_{1}%
^{1+h}\left(  \sigma\left(  Y_{s}^{1,x,c}\right)  \right)  ^{2}ds\right]
\right)  ^{\frac{1}{2}}\\
\leq &  C_{0}\sqrt{\underline{\sigma}^{2}+\overline{\sigma}^{2}}\sqrt
{h}\text{.\ \ \ \ \ \ }\blacksquare
\end{align*}
\hfill\hfill}

\begin{lemma}
{\normalsize \label{lemma-taylor-O} Let $\{H_{t}\}_{t\in\lbrack0,1]}$ be the
functions defined in (\ref{hf}),
%{\small [}$\{X_{i}\}\,${\small \ a sequence
%of random variables with conditional upper and lower variances }%
%$\overline{\sigma}^{2}${\small  and }$\underline{\sigma}^{2}${\small
%respectively,]}
and define the family of functions $\{L_{m,n}\}_{m=1}^{n}$ by
%The set P, (and upper expectation), enters here, while it did not in Lemma A.1%
\begin{equation}
L_{m,n}(x)=H_{\frac{m}{n}}(x)+\frac{\overline{\sigma}^{2}}{2n}\left[
{H_{\frac{m}{n}}^{\prime\prime}(x)}\right]  ^{+}-\frac{\underline{\sigma}^{2}%
}{2n}\left[  {H_{\frac{m}{n}}^{\prime\prime}(x)}\right]  ^{-}. \label{Lt}%
\end{equation}
{ }Then
\begin{equation}
\lim_{n\rightarrow\infty}\sum_{m=1}^{n}\left\vert \mathbb{E}\left[
H_{\frac{m}{n}}\left(  \frac{\Sigma_{1}^{m}X_{i}}{\sqrt{n}}\right)  \right]
-\mathbb{E}\left[  L_{m,n}\left(  \frac{\Sigma_{1}^{m-1}X_{i}}{\sqrt{n}%
}\right)  \right]  \right\vert =0. \label{lemma-taylor1}%
\end{equation}
}
\end{lemma}

{\normalsize \noindent\textbf{Proof:} It suffices to prove
\begin{align}
\lim\limits_{n\rightarrow\infty}\sum_{m=1}^{n}\left\vert \mathbb{E}\left[
H_{\frac{m}{n}}\left(  \frac{\Sigma_{1}^{m}X_{i}}{\sqrt{n}}\right)  \right]
-f(m,n)\right\vert  &  =0\text{ and}\label{O-var-taylor1}\\
\lim\limits_{n\rightarrow\infty}\sum_{m=1}^{n}\left\vert f(m,n)-\mathbb{E}%
\left[  L_{m,n}\left(  \frac{\Sigma_{1}^{m-1}X_{i}}{\sqrt{n}}\right)  \right]
\right\vert  &  =0, \label{O-var-taylor2}%
\end{align}
where
\[
f(m,n)=\mathbb{E}\left[  H_{\frac{m}{n}}\left(  \frac{\Sigma_{1}^{m-1}X_{i}%
}{\sqrt{n}}\right)  +H_{\frac{m}{n}}^{\prime}\left(  \frac{\Sigma_{1}%
^{m-1}X_{i}}{\sqrt{n}}\right)  \frac{X_{m}}{\sqrt{n}}+H_{\frac{m}{n}}%
^{\prime\prime}\left(  \frac{\Sigma_{1}^{m-1}X_{i}}{\sqrt{n}}\right)
\frac{X_{m}^{2}}{2n}\right]  .
\]
}

{\normalsize By Lemma \ref{lemma-ddp}, there exists $L>0$ such that
\[
\sup\limits_{t\in\lbrack0,1]}\sup\limits_{x\in\mathbb{R}}\left\vert
H_{t}^{\prime\prime}(x)\right\vert \leq L,\ \sup\limits_{t\in\lbrack0,1]}%
\sup\limits_{x,y\in\mathbb{R},x\neq y}\frac{\left\vert H_{t}^{\prime\prime
}(x)-H_{t}^{\prime\prime}(y)\right\vert }{|x-y|}\leq L.
\]
By the Taylor expansion of $H_{t}\in C_{b}^{2}(\mathbb{R})$, $\forall
\epsilon>0$ $\exists\delta>0$ ($\delta$ depends only on $L$ and $\epsilon$),
such that, for any $x,y\in\mathbb{R}$ and $t\in\lbrack0,1]$,
\begin{equation}
\left\vert H_{t}(x+y)-H_{t}(x)-H_{t}^{\prime}(x)y-\frac{1}{2}H_{t}%
^{\prime\prime}(x)y^{2}\right\vert \leq\epsilon|y|^{2}I_{\{|y|<\delta
\}}+L|y|^{2}I_{\{|y|\geq\delta\}}. \label{le0}%
\end{equation}
Let $x=\Sigma_{1}^{m-1}X_{i}/\sqrt{n}$ and $y=X_{m}/\sqrt{n}$ in (\ref{le0})
to derive
\[
\sum_{m=1}^{n}\left\vert \mathbb{E}\left[  H_{\frac{m}{n}}\left(  \frac
{\Sigma_{1}^{m}X_{i}}{\sqrt{n}}\right)  \right]  -f(m,n)\right\vert
\leq\overline{\sigma}^{2}\epsilon+\frac{L}{n}\sum_{m=1}^{n}\mathbb{E}\left[
|X_{m}|^{2}I_{\{|X_{m}|\geq\sqrt{n}\delta\}}\right]  .
\]
By the arbitrariness of $\epsilon$ and the Lindeberg condition (\ref{linder}),
we obtain (\ref{O-var-taylor1}). }

{\normalsize By Lemma \ref{lemma-rect}, we have%
\begin{align*}
f(m,n)=  &  \mathbb{E}\left[  H_{\frac{m}{n}}\left(  \frac{\Sigma_{1}%
^{m-1}X_{i}}{\sqrt{n}}\right)  +H_{\frac{m}{n}}^{\prime}\left(  \frac
{\Sigma_{1}^{m-1}X_{i}}{\sqrt{n}}\right)  \frac{X_{m}}{\sqrt{n}}+H_{\frac
{m}{n}}^{\prime\prime}\left(  \frac{\Sigma_{1}^{m-1}X_{i}}{\sqrt{n}}\right)
\frac{X_{m}^{2}}{2n}\right] \\
=  &  \mathbb{E}\left[  H_{\frac{m}{n}}\left(  \frac{\Sigma_{1}^{m-1}X_{i}%
}{\sqrt{n}}\right)  +H_{\frac{m}{n}}^{\prime\prime}\left(  \frac{\Sigma
_{1}^{m-1}X_{i}}{\sqrt{n}}\right)  \frac{X_{m}^{2}}{2n}\right] \\
=  &  \mathbb{E}\left[  H_{\frac{m}{n}}\left(  \frac{\Sigma_{1}^{m-1}X_{i}%
}{\sqrt{n}}\right)  +\frac{1}{2n}\mathbb{E}\left[  H_{\frac{m}{n}}%
^{\prime\prime}\left(  \frac{\Sigma_{1}^{m-1}X_{i}}{\sqrt{n}}\right)
X_{m}^{2}|\mathcal{G}_{m-1}\right]  \right] \\
=  &  \mathbb{E}\left[  L_{m,n}\left(  \frac{\Sigma_{1}^{m-1}X_{i}}{\sqrt{n}%
}\right)  \right]  .
\end{align*}
\noindent This implies (\ref{O-var-taylor2}) and completes the proof.
\hfill$\blacksquare$ \medskip}

\subsection{Proof of the CLT (Theorem \ref{thm-CLT})\label{app-CLTproof}}

{\normalsize {For $h>0$ sufficiently small, }let $\{H_{t}\}_{t\in
\lbrack0,1+h]}$ be the functions defined by (\ref{hf}). First prove%
\[
\lim_{n\rightarrow\infty}\left\vert \mathbb{E}\left[  H_{1}\left(
\frac{\Sigma_{1}^{n}X_{i}}{\sqrt{n}}\right)  \right]  -E_{P^{\ast}}\left[
\varphi\left(  W_{1+h}^{c}\right)  \right]  \right\vert =0\text{.}%
\]
We have%
\begin{align*}
&  \mathbb{E}\left[  H_{1}\left(  \frac{\Sigma_{1}^{n}X_{i}}{\sqrt{n}}\right)
\right]  -E_{P^{\ast}}\left[  \varphi\left(  W_{1+h}^{c}\right)  \right] \\
=  &  \mathbb{E}\left[  H_{1}\left(  \frac{\Sigma_{1}^{n}X_{i}}{\sqrt{n}%
}\right)  \right]  -H_{0}(0)\\
=  &  \sum\limits_{m=1}^{n}\left\{  \mathbb{E}\left[  H_{\frac{m}{n}}\left(
\frac{\Sigma_{1}^{m}X_{i}}{\sqrt{n}}\right)  \right]  -\mathbb{E}\left[
H_{\frac{m-1}{n}}\left(  \frac{\Sigma_{1}^{m-1}X_{i}}{\sqrt{n}}\right)
\right]  \right\} \\
=  &  \sum\limits_{m=1}^{n}\left\{  \mathbb{E}\left[  H_{\frac{m}{n}}\left(
\frac{\Sigma_{1}^{m}X_{i}}{\sqrt{n}}\right)  \right]  -\mathbb{E}\left[
L_{m,n}\left(  \frac{\Sigma_{1}^{m-1}X_{i}}{\sqrt{n}}\right)  \right]
\right\} \\
&  +\sum_{m=1}^{n}\left\{  \mathbb{E}\left[  L_{m,n}\left(  \frac{\Sigma
_{1}^{m-1}X_{i}}{\sqrt{n}}\right)  \right]  -\mathbb{E}\left[  H_{\frac
{m-1}{n}}\left(  \frac{\Sigma_{1}^{m-1}X_{i}}{\sqrt{n}}\right)  \right]
\right\} \\
=  &  :I_{1n}+I_{2n}\text{,}%
\end{align*}
where $L_{m,n}(x)=H_{\frac{m}{n}}(x)+\frac{\overline{\sigma}^{2}}{2n}\left[
{H_{\frac{m}{n}}^{\prime\prime}(x)}\right]  ^{+}-\frac{\underline{\sigma}^{2}%
}{2n}\left[  {H_{\frac{m}{n}}^{\prime\prime}(x)}\right]  ^{-},\ 1\leq m\leq
n$. }

{\normalsize By Lemma \ref{lemma-taylor-O}, \
\[
|I_{1n}|\leq\sum\limits_{m=1}^{n}\left\vert \mathbb{E}\left[  H_{\frac{m}{n}%
}\left(  \frac{\Sigma_{1}^{m}X_{i}}{\sqrt{n}}\right)  \right]  -\mathbb{E}%
\left[  L_{m,n}\left(  \frac{\Sigma_{1}^{m-1}X_{i}}{\sqrt{n}}\right)  \right]
\right\vert \rightarrow0\text{ ~as }n\rightarrow\infty\text{.}%
\]
Furthermore, by Lemma \ref{lemma-ddp}(5), as $n\rightarrow\infty$,
\begin{align*}
|I_{2n}|  &  \leq\sum_{m=1}^{n}\mathbb{E}\left[  \left\vert L_{m,n}\left(
\frac{\Sigma_{1}^{m-1}X_{i}}{\sqrt{n}}\right)  -H_{\frac{m-1}{n}}\left(
\frac{\Sigma_{1}^{m-1}X_{i}}{\sqrt{n}}\right)  \right\vert \right] \\
&  \leq\sum_{m=1}^{n}\sup\limits_{x\in\mathbb{R}}\left\vert L_{m,n}%
(x)-H_{\frac{m-1}{n}}(x)\right\vert \\
&  =\sum_{m=1}^{n}\sup\limits_{x\in\mathbb{R}}\left\vert H_{\frac{m-1}{n}%
}(x)-H_{\frac{m}{n}}\left(  x\right)  -\frac{\overline{\sigma}^{2}}{2n}\left[
{H_{\frac{m}{n}}^{\prime\prime}(x)}\right]  ^{+}+\frac{\underline{\sigma}^{2}%
}{2n}\left[  {H_{\frac{m}{n}}^{\prime\prime}(x)}\right]  ^{-}\right\vert
\rightarrow0\text{.}%
\end{align*}
By Lemma \ref{lemma-ddp}(6),
%$$\sup_{x\in\mathbb{R}}\left|H_1(x)-\varphi(x)\right|\le C_0\sqrt{\overline{\sigma}^2+\underline{\sigma}^2}\sqrt{h}. $$
%Similarly we have,
%\begin{align*}
%&\left|E_{P^*}\left[  \varphi\left(  W_{1+h}^{c}\right)  \right]-E_{P^*}\left[  \varphi\left(  W_{1}^{c}\right)  \right]\right|\\
%\le  &\|\varphi'\|E_{P^*}\left[  \left|  W_{1+h}^{c}-  W_{1}^{c} \right|\right]\\
%\le  &\|\varphi'\|E_{P^*}\left[  \left| \int_1^{1+h} \sigma \left(W_{s}^{c}\right)dB_s \right|\right]\\
%\le& \|\varphi'\|\sqrt{\overline{\sigma}^2+\underline{\sigma}^2}\sqrt{h},\
%\end{align*}
%Then we have%
$\lim_{n\rightarrow\infty}\left\vert \mathbb{E}\left[  \varphi\left(
\frac{\Sigma_{1}^{n}X_{i}}{\sqrt{n}}\right)  \right]  -E_{P^{\ast}}\left[
\varphi\left(  W_{1}^{c}\right)  \right]  \right\vert \leq$\newline%
$\lim_{n\rightarrow\infty}\left\vert \mathbb{E}\left[  \varphi\left(
\frac{\Sigma_{1}^{n}X_{i}}{\sqrt{n}}\right)  \right]  -\mathbb{E}\left[
H_{1}\left(  \frac{\Sigma_{1}^{n}X_{i}}{\sqrt{n}}\right)  \right]  \right\vert
+$\newline$\lim_{n\rightarrow\infty}\left\vert \mathbb{E}\left[  H_{1}\left(
\frac{\Sigma_{1}^{n}X_{i}}{\sqrt{n}}\right)  \right]  -E_{P^{\ast}}\left[
\varphi\left(  W_{1+h}^{c}\right)  \right]  \right\vert +\left\vert
E_{P^{\ast}}\left[  \varphi\left(  W_{1+h}^{c}\right)  \right]  -E_{P^{\ast}%
}\left[  \varphi\left(  W_{1}^{c}\right)  \right]  \right\vert $\newline%
$\leq\sup_{x\in\mathbb{R}}\left\vert H_{1}\left(  x\right)  -\varphi\left(
x\right)  \right\vert +C_{0}\sqrt{\overline{\sigma}^{2}+\underline{\sigma}%
^{2}}\sqrt{h}~\leq~2C_{0}\sqrt{\overline{\sigma}^{2}+\underline{\sigma}^{2}%
}\sqrt{h}$.\newline Since $h$ is arbitrary, the proof is complete.
\medskip\hfill$\blacksquare$ \medskip}

{\normalsize
%\begin{lemma}\label{ph-clt-exist}
%\begin{description}
%\item[(1).] For any $\varphi\in\Psi(\underline{\sigma},\overline{\sigma},c)$,  which satisfies condition \textup{(C1)} or \textup{(C2)}, there exists a sequence of probability measures   $\{Q_{n}^c: Q_{n}^c\in\mathcal{P},n\geq1\}$ such that
%\begin{equation}\label{ph-exist-clt1}
%\lim\limits_{n\rightarrow\infty}E_{Q_n^c}\left[
%\varphi \left(\frac{S_n}{\sqrt{n}} \right)\right]
%=E_{P^*}\left[\varphi\left(W_1^{\underline{\sigma},\overline{\sigma},c}\right)\right],
%\end{equation}
%\item[(2).] For any $\varphi\in\Psi(\overline{\sigma},\underline{\sigma},c)$,  which satisfies condition \textup{(C1)} or \textup{(C2)}, there exists a sequence of probability measures   $\{\Q_{n}^c: \Q_{n}^c\in\mathcal{P},n\geq1\}$ such that
%\begin{equation}\label{ph-exist-clt2}
%\lim\limits_{n\rightarrow\infty}E_{\Q_n^c}\left[
%\varphi \left(\frac{S_n}{\sqrt{n}} \right)\right]
%=E_{P^*}\left[\varphi\left(W_1^{c}\right)\right],
%\end{equation}
%\end{description}
%\end{lemma}
\
%\begin{proof}[\bf Proof]
%With  Lemma \ref{lemma-taylor}-(2)(3), the proof method in Theorem \ref{thm-CLT-G} is also applicable here, so we omit it.
%\end{proof}
\noindent}

\subsection{A corollary\label{app-indicators}}

Indicator functions for one-sided intervals $[c,\infty)$ can be suitably
approximated by functions $\varphi$ satisfying the conditions in Theorem
\ref{thm-CLT}, which suggests that the limiting result
(\ref{thm-clt-eq-O-leq-sup}) is valid also for such indicators. The following
corollary confirms this, and is of interest also because it is used below in
the proof of Theorem \ref{thm-nolearn}. See our working paper version
(Corollary 3.4) for a more general result that considers also indicators for
intervals of the form $(-\infty,c]$.

\begin{corol}
\label{cor-pro-O}Adopt the assumptions in Theorem \ref{thm-CLT}. Then, for any
$c\in\mathbb{R}$,
\begin{equation}%
\begin{array}
[c]{ll}%
\lim\limits_{n\rightarrow\infty}\sup_{Q\in\mathcal{P}}Q\left(  \frac
{\Sigma_{1}^{n}X_{i}}{\sqrt{n}}\geq c\right)  & =P^{\ast}\left(  W_{1}^{c}\geq
c\right)
\end{array}
\label{geq-clt}%
\end{equation}
and%
\begin{equation}
P^{\ast}\left(  W_{1}^{c}\geq c\right)  =\left\{
\begin{array}
[c]{ccc}%
\frac{2\overline{\sigma}}{\overline{\sigma}+\underline{\sigma}}\Phi\left(
-\frac{c}{\overline{\sigma}}\right)  &  & c>0\\
1-\frac{2\underline{\sigma}}{\overline{\sigma}+\underline{\sigma}}\Phi\left(
\frac{c}{\underline{\sigma}}\right)  &  & c\leq0
\end{array}
\right.  \text{,} \label{Pindicator}%
\end{equation}
where $\Phi$ is the standard normal cdf.
\end{corol}

\noindent\textbf{Proof: }{\normalsize {For any $c\in\mathbb{R}$ and
$\varepsilon>0$, suppose that $f_{1},g_{1}\in C_{b}^{3}(\mathbb{R}_{+})$
satisfy
\[
\left\{
\begin{array}
[c]{l}%
f_{1}(x)=1\text{ \ for }x\geq\frac{\underline{\sigma}}{\overline{\sigma}%
}\varepsilon\\
f_{1}^{\prime\prime}(x)\leq0\text{ \ for }x\geq0\\
f_{1}(0)=\frac{\overline{\sigma}}{\overline{\sigma}+\underline{\sigma}}%
\end{array}
\right.  \left\{
\begin{array}
[c]{ll}%
g_{1}(x)=1\text{ \ for }x\geq\varepsilon & \\
g_{1}^{\prime\prime}(x)\leq0\text{ \ for }x\geq0 & \\
g_{1}(0)=\frac{\overline{\sigma}}{\overline{\sigma}+\underline{\sigma}} &
\end{array}
\right.
\]
Define $f_{\varepsilon}$ and $g_{\varepsilon}$ by
\begin{equation}
f_{\varepsilon}(x)=\left\{
\begin{array}
[c]{ll}%
f_{1}(x-c-\varepsilon) & \text{for }x\geq c+\varepsilon\\
-\frac{\overline{\sigma}}{\underline{\sigma}}f_{1}\left(  -\frac
{\underline{\sigma}}{\overline{\sigma}}(x-c-\varepsilon)\right)
+\frac{\overline{\sigma}}{\underline{\sigma}} & \text{for }x\leq c+\varepsilon
\end{array}
\right.  \label{f-ep}%
\end{equation}%
\begin{equation}
g_{\varepsilon}(x)=\left\{
\begin{array}
[c]{ll}%
g_{1}(x-c+\varepsilon) & \text{for }x\geq c-\varepsilon\\
-\frac{\overline{\sigma}}{\underline{\sigma}}g_{1}\left(  -\frac
{\underline{\sigma}}{\overline{\sigma}}(x-c+\varepsilon)\right)
+\frac{\overline{\sigma}}{\underline{\sigma}} & \text{for }x\leq c-\varepsilon
\end{array}
\right.  \label{g-ep}%
\end{equation}
} It can be checked that
\begin{align*}
g_{\varepsilon}(x)  &  \geq I_{[c,\infty)}(x)\geq f_{\varepsilon}(x)\text{
\ and }\\
|g_{\varepsilon}(x)-f_{\varepsilon}(x)|  &  \leq I_{\left[  c-(1+\frac
{\overline{\sigma}}{\underline{\sigma}})\varepsilon,c+(1+\frac{\underline
{\sigma}}{\overline{\sigma}})\varepsilon\right]  }(x),\ \forall x\in
\mathbb{R}.
\end{align*}
Consider the solution $(\widetilde{W}_{t}^{x})_{t\geq0}$ of the SDE
\begin{equation}
\left\{
\begin{array}
[c]{l}%
d\widetilde{W}_{t}^{x}=\left(  \underline{\sigma}I_{[0,\infty)}\left(
\widetilde{W}_{t}^{x}\right)  +\overline{\sigma}I_{(-\infty,0)}\left(
\widetilde{W}_{t}^{x}\right)  \right)  dB_{t},\quad t\geq0\\
\widetilde{W}_{0}^{x}=x.
\end{array}
\right.  \label{sde4}%
\end{equation}
Then $W_{1}^{c}$ and $c+\widetilde{W}_{1}^{-c}$ are described by the same law,
and }\newline{\small
\begin{align*}
&  \left\vert \sup\limits_{Q\in\mathcal{P}}Q\left(  \frac{\Sigma_{1}^{n}X_{i}%
}{\sqrt{n}}\geq c\right)  -P^{\ast}\left(  W_{1}^{c}\geq c\right)  \right\vert
\\
\leq &  \left\vert \mathbb{E}\left[  f_{\varepsilon}\left(  \frac{\Sigma
_{1}^{n}X_{i}}{\sqrt{n}}\right)  \right]  -E_{P^{\ast}}\left[  g_{\varepsilon
}\left(  W_{1}^{c}\right)  \right]  \right\vert +\left\vert \mathbb{E}\left[
g_{\varepsilon}\left(  \frac{\Sigma_{1}^{n}X_{i}}{\sqrt{n}}\right)  \right]
-E_{P^{\ast}}\left[  f_{\varepsilon}\left(  W_{1}^{c}\right)  \right]
\right\vert \\
\leq &  \left\vert \mathbb{E}\left[  f_{\varepsilon}\left(  \frac{\Sigma
_{1}^{n}X_{i}}{\sqrt{n}}\right)  \right]  -E_{P^{\ast}}\left[  f_{\varepsilon
}\left(  W_{1}^{c+\varepsilon}\right)  \right]  \right\vert +\left\vert
\mathbb{E}\left[  g_{\varepsilon}\left(  \frac{\Sigma_{1}^{n}X_{i}}{\sqrt{n}%
}\right)  \right]  -E_{P^{\ast}}\left[  g_{\varepsilon}\left(  W_{1}%
^{c-\varepsilon}\right)  \right]  \right\vert \\
&  +\left\vert E_{P^{\ast}}\left[  f_{\varepsilon}\left(  W_{1}^{c+\varepsilon
}\right)  -f_{\varepsilon}\left(  W_{1}^{c}\right)  \right]  \right\vert
+\left\vert E_{P^{\ast}}\left[  g_{\varepsilon}\left(  W_{1}^{c-\varepsilon
}\right)  -g_{\varepsilon}\left(  W_{1}^{c}\right)  \right]  \right\vert \\
&  +2\left\vert E_{P^{\ast}}\left[  f_{\varepsilon}\left(  W_{1}^{c}\right)
-g_{\varepsilon}\left(  W_{1}^{c}\right)  \right]  \right\vert \\
\leq &  \left\vert \mathbb{E}\left[  f_{\varepsilon}\left(  \frac{\Sigma
_{1}^{n}X_{i}}{\sqrt{n}}\right)  \right]  -E_{P^{\ast}}\left[  f_{\varepsilon
}\left(  W_{1}^{c+\varepsilon}\right)  \right]  \right\vert +\left\vert
\mathbb{E}\left[  g_{\varepsilon}\left(  \frac{\Sigma_{1}^{n}X_{i}}{\sqrt{n}%
}\right)  \right]  -E_{P^{\ast}}\left[  g_{\varepsilon}\left(  W_{1}%
^{c-\varepsilon}\right)  \right]  \right\vert \\
&  +E_{P^{\ast}}\left[  \left\vert f_{\varepsilon}\left(  c+\varepsilon
+\widetilde{W}_{1}^{-c-\varepsilon}\right)  -f_{\varepsilon}\left(
c+\widetilde{W}_{1}^{-c}\right)  \right\vert +\left\vert g_{\varepsilon
}\left(  c-\varepsilon+\widetilde{W}_{1}^{-c+\varepsilon}\right)
-g_{\varepsilon}\left(  c+\widetilde{W}_{1}^{-c}\right)  \right\vert \right]
\\
&  +2\left\vert E_{P^{\ast}}\left[  f_{\varepsilon}\left(  W_{1}^{c}\right)
-g_{\varepsilon}\left(  W_{1}^{c}\right)  \right]  \right\vert \\
\leq &  \left\vert \mathbb{E}\left[  f_{\varepsilon}\left(  \frac{\Sigma
_{1}^{n}X_{i}}{\sqrt{n}}\right)  \right]  -E_{P^{\ast}}\left[  f_{\varepsilon
}\left(  W_{1}^{c+\varepsilon}\right)  \right]  \right\vert +\left\vert
\mathbb{E}\left[  g_{\varepsilon}\left(  \frac{\Sigma_{1}^{n}X_{i}}{\sqrt{n}%
}\right)  \right]  -E_{P^{\ast}}\left[  g_{\varepsilon}\left(  W_{1}%
^{c-\varepsilon}\right)  \right]  \right\vert \\
&  +C_{0}E_{P^{\ast}}\left[  2\varepsilon+\left\vert \widetilde{W}%
_{1}^{-c-\varepsilon}-\widetilde{W}_{1}^{-c}\right\vert +\left\vert
\widetilde{W}_{1}^{-c+\varepsilon}-\widetilde{W}_{1}^{-c}\right\vert \right]
+2P^{\ast}\left(  c-\left(  1+\frac{\overline{\sigma}}{\underline{\sigma}%
}\right)  \varepsilon\leq W_{1}^{c}\leq c+\left(  1+\frac{\underline{\sigma}%
}{\overline{\sigma}}\right)  \varepsilon\right)  \text{,}%
\end{align*}
}{\normalsize where $C_{0}$ is a constant that depends on $\Vert
f_{\varepsilon}^{\prime}\Vert,\Vert g_{\varepsilon}^{\prime}\Vert$. With Le
Gall (1984, Theorem 1.5) and Theorem \ref{thm-CLT}, the upper probability
equation in (\ref{geq-clt}) is proven. \medskip}

The expression (\ref{Pindicator}) may be derived by integrating the pdf
(\ref{pdf-c>0})-(\ref{pdf-c<0}). {\normalsize \hfill$\blacksquare$ \medskip}

\section{Appendix: Proofs for bandits\label{app-bandits}}

\subsection{An explicit density\label{app-density}}

Let $W_{1}^{c}$ be the $t=1$ value of the oscillating Brownian motion defined
by (\ref{sde0})-(\ref{sigma-function}). Keilson and Wellner (1978, Theorem 1)
give the following expression for its pdf: For $c\geq0$,%

\begin{equation}
q\left(  y\right)  =\left\{
\begin{array}
[c]{ccc}%
\frac{1}{\underline{\sigma}\sqrt{2\pi}}\left(  e^{-\frac{(\frac{-c}%
{\overline{\sigma}}-\frac{y-c}{\underline{\sigma}})^{2}}{2}}+\frac
{\overline{\sigma}-\underline{\sigma}}{\underline{\sigma}+\overline{\sigma}%
}e^{-\frac{(\frac{c}{\overline{\sigma}}+\frac{y-c}{\underline{\sigma}})^{2}%
}{2}}\right)  &  & y\geq c\\
&  & \\
\frac{1}{\overline{\sigma}\sqrt{2\pi}}\left(  e^{-\frac{(\frac{-c}%
{\overline{\sigma}}-\frac{y-c}{\overline{\sigma}})^{2}}{2}}-\frac
{\overline{\sigma}-\underline{\sigma}}{\underline{\sigma}+\overline{\sigma}%
}e^{-\frac{(\frac{c}{\overline{\sigma}}+\frac{c-y}{\overline{\sigma}})^{2}}%
{2}}\right)  &  & y<c
\end{array}
\right.  \label{pdf-c>0}%
\end{equation}
and for $c<0$,%

\begin{equation}
q\left(  y\right)  =\left\{
\begin{array}
[c]{ccc}%
\frac{1}{\underline{\sigma}\sqrt{2\pi}}\left(  e^{-\frac{(\frac{-c}%
{\underline{\sigma}}-\frac{y-c}{\underline{\sigma}})^{2}}{2}}+\frac
{\overline{\sigma}-\underline{\sigma}}{\underline{\sigma}+\overline{\sigma}%
}e^{-\frac{(\frac{-c}{\underline{\sigma}}+\frac{y-c}{\underline{\sigma}})^{2}%
}{2}}\right)  &  & y\geq c\\
&  & \\
\frac{1}{\overline{\sigma}\sqrt{2\pi}}\left(  e^{-\frac{(\frac{-c}%
{\overline{\sigma}}-\frac{y-c}{\overline{\sigma}})^{2}}{2}}-\frac
{\overline{\sigma}-\underline{\sigma}}{\underline{\sigma}+\overline{\sigma}%
}e^{-\frac{(\frac{c}{\overline{\sigma}}+\frac{c-y}{\overline{\sigma}})^{2}}%
{2}}\right)  &  & y<c
\end{array}
\right.  \label{pdf-c<0}%
\end{equation}

\noindent These expressions are used to derive (\ref{Vexp}) and to prove
Corollary \ref{cor-pro-O} and Theorem \ref{thm-V}.

\subsection{\noindent Proof of Theorem \ref{thm-V}}

As indicated in the text, (i) follows from Theorem \ref{thm-CLT} and the above
density; and (ii) follows from (i) by inspection of the above density. It
remains to prove (iii).

\smallskip Take $c\geq0$. The proof for $c<0$ is similar. In light of
(\ref{equiv}) and (\ref{thm-clt-eq-O-leq-sup}), it suffices to compute
$E_{P^{\ast}}[\varphi(W_{1}^{c})]$. Use the pdf of $W_{1}^{c}$ in
(\ref{pdf-c>0}), to deduce that, for $c\geq0$,%
\begin{align*}
E_{P^{\ast}}\left[  \varphi\left(  W_{1}^{c}\right)  \right]   &  =\int
_{c}^{\infty}q\left(  y\right)  \varphi_{1}\left(  y-c\right)  dy+\int
\limits_{-\infty}^{c}q\left(  y\right)  \left[  -\frac{\overline{\sigma}%
}{\underline{\sigma}}\varphi_{1}\left(  -\frac{\underline{\sigma}}%
{\overline{\sigma}}(y-c)\right)  \right]  dy\\
&  =\frac{1}{\sqrt{2\pi}}\frac{1}{\underline{\sigma}}\int\limits_{c}^{\infty
}\varphi_{1}\left(  y-c\right)  \left[  \frac{2\overline{\sigma}}%
{\overline{\sigma}+\underline{\sigma}}\right]  e^{-\frac{\left(  \frac
{c}{\overline{\sigma}}+\frac{y-c}{\underline{\sigma}}\right)  ^{2}}{2}}dy\\
&  -\frac{1}{\sqrt{2\pi}}\frac{1}{\underline{\sigma}}\int\limits_{-\infty}%
^{c}\varphi_{1}\left(  -\frac{\underline{\sigma}}{\overline{\sigma}%
}(y-c)\right)  \left(  e^{-\frac{(\frac{-c}{\overline{\sigma}}-\frac
{y-c}{\overline{\sigma}})^{2}}{2}}-\frac{\overline{\sigma}-\underline{\sigma}%
}{\underline{\sigma}+\overline{\sigma}}e^{-\frac{(\frac{c}{\overline{\sigma}%
}+\frac{c-y}{\overline{\sigma}})^{2}}{2}}\right)  dy\ \ \\
&  =\frac{1}{\sqrt{2\pi}}\frac{1}{\underline{\sigma}}\int\limits_{0}^{\infty
}\varphi_{1}\left(  y\right)  \left[  \frac{2\overline{\sigma}}{\overline
{\sigma}+\underline{\sigma}}\right]  e^{-\frac{\left(  \frac{c}{\overline
{\sigma}}+\frac{y}{\underline{\sigma}}\right)  ^{2}}{2}}dy\\
&  -\frac{1}{\sqrt{2\pi}}\frac{1}{\underline{\sigma}}\int\limits_{-\infty}%
^{0}\varphi_{1}\left(  -\frac{\underline{\sigma}}{\overline{\sigma}}y\right)
\left(  e^{-\frac{(\frac{-c}{\overline{\sigma}}-\frac{y}{\overline{\sigma}%
})^{2}}{2}}-\frac{\overline{\sigma}-\underline{\sigma}}{\underline{\sigma
}+\overline{\sigma}}e^{-\frac{(\frac{c}{\overline{\sigma}}-\frac{y}%
{\overline{\sigma}})^{2}}{2}}\right)  dy\ \ \\
&  =\frac{1}{\sqrt{2\pi}}\frac{1}{\underline{\sigma}}\int\limits_{0}^{\infty
}\varphi_{1}\left(  y\right)  \left[  \frac{2\overline{\sigma}}{\overline
{\sigma}+\underline{\sigma}}\right]  e^{-\frac{\left(  \frac{c}{\overline
{\sigma}}+\frac{y}{\underline{\sigma}}\right)  ^{2}}{2}}dy\\
&  -\frac{1}{\sqrt{2\pi}}\frac{1}{\underline{\sigma}}\frac{\overline{\sigma}%
}{\underline{\sigma}}\int\limits_{0}^{\infty}\varphi_{1}\left(  z\right)
\left(  e^{-\frac{(\frac{-c}{\overline{\sigma}}+\frac{z}{\overline
{\underline{\sigma}}})^{2}}{2}}-\frac{\overline{\sigma}-\underline{\sigma}%
}{\underline{\sigma}+\overline{\sigma}}e^{-\frac{(\frac{c}{\overline{\sigma}%
}+\frac{z}{\underline{\sigma}})^{2}}{2}}\right)  dz\text{ }\\
&  =\frac{1}{\sqrt{2\pi}\underline{\sigma}}\int_{0}^{\infty}\varphi_{1}\left(
y\right)  \frac{\overline{\sigma}}{\underline{\sigma}}\left[  e^{-\frac
{\left(  y+m\right)  ^{2}}{2\underline{\sigma}^{2}}}-e^{-\frac{(y-m)^{2}%
}{2\underline{\sigma}^{2}}}\right]  dy
\end{align*}
where $m=\frac{\underline{\sigma}}{\overline{\sigma}}c$. Thus we want to prove
that
\[
\frac{1}{\sqrt{2\pi}\underline{\sigma}}\int_{0}^{\infty}\varphi_{1}\left(
y\right)  \frac{\overline{\sigma}}{\underline{\sigma}}\left[  e^{-\frac
{\left(  y+m\right)  ^{2}}{2\underline{\sigma}^{2}}}-e^{-\frac{(y-m)^{2}%
}{2\underline{\sigma}^{2}}}\right]  dy\geq-\frac{\overline{\sigma}}%
{\underline{\sigma}}\varphi_{1}\left(  \frac{\underline{\sigma}}%
{\overline{\sigma}}c\right)  \text{,}%
\]

\noindent with equality if and only if $c=0$. \noindent\ $.$

It is evident that $E_{P^{\ast}}\left[  \varphi\left(  W_{1}^{c}\right)
\right]  =0=\varphi\left(  0\right)  $ if $c=0$. Henceforth, take $c>0$ and
prove that
\[
\int_{0}^{\infty}\varphi_{1}\left(  y\right)  \left[  \left(  e^{-\frac
{\left(  y-m\right)  ^{2}}{2\underline{\sigma}^{2}}}-e^{-\frac{(y+m)^{2}%
}{2\underline{\sigma}^{2}}}\right)  /\sqrt{2\pi}\underline{\sigma}\right]
dy<~\varphi_{1}\left(  m\right)  \text{.}%
\]

Denote by $f\left(  y\right)  $ the expression in the square bracket, (thus
$f\left(  y\right)  >0$ for all $y>0$), and let $F\equiv%
%TCIMACRO{\dint \limits_{0}^{\infty}}%
%BeginExpansion
{\displaystyle\int\limits_{0}^{\infty}}
%EndExpansion
f\left(  y\right)  dy$, $0<F<1$. Then $f/F$ is a density. If its mean is $\mu
$, then, by strict concavity of $\varphi_{1}$,
\begin{equation}
\int_{0}^{\infty}\varphi_{1}\left(  y\right)  f\left(  y\right)
dy<F\varphi_{1}\left(  \mu\right)  \text{.} \label{V1}%
\end{equation}

\noindent Next we prove that $F\mu=m$:
\begin{align*}
F\mu &  =\int_{0}^{\infty}y\left[  \left(  e^{-\frac{\left(  y-m\right)  ^{2}%
}{2\underline{\sigma}^{2}}}-e^{-\frac{(y+m)^{2}}{2\underline{\sigma}^{2}}%
}\right)  /\sqrt{2\pi}\underline{\sigma}\right]  dy\\
&  =\int_{-m}^{\infty}(z+m)\left[  e^{-\frac{z^{2}}{2\underline{\sigma}^{2}}%
}/\sqrt{2\pi}\underline{\sigma}\right]  dz-\int_{m}^{\infty}(z-m)\left[
\left(  e^{-\frac{z^{2}}{2\underline{\sigma}^{2}}}\right)  /\sqrt{2\pi
}\underline{\sigma}\right]  dz\newline\\
&  =\int_{-m}^{m}z\left[  e^{-\frac{z^{2}}{2\underline{\sigma}^{2}}}%
/\sqrt{2\pi}\underline{\sigma}\right]  dz+m\int_{-m}^{\infty}\left[
e^{-\frac{z^{2}}{2\underline{\sigma}^{2}}}/\sqrt{2\pi}\underline{\sigma
}\right]  dz+m\int_{m}^{\infty}\left[  e^{-\frac{z^{2}}{2\underline{\sigma
}^{2}}}/\sqrt{2\pi}\underline{\sigma}\right]  dz\\
&  =0+m\left[  \Pr\left(  Z>-m\right)  +\Pr\left(  Z>m\right)  \right]
~=m\text{,}%
\end{align*}
where probabilities are computed according to $\mathbb{N}\left(
0,\underline{\sigma}^{2}\right)  $.

Finally, $F\mu=m$ $~\Longrightarrow$ $F\varphi_{1}\left(  \mu\right)
=F\varphi_{1}\left(  m/F\right)  \leq\varphi_{1}\left(  m\right)  $, by $F<1$,
$\varphi_{1}\left(  0\right)  =0$, and the concavity of $\varphi_{1}$. Combine
with (\ref{V1}) to complete the proof.{\normalsize \hfill$\blacksquare$
\medskip}

\subsection{\noindent Proof of Theorem \ref{thm-nolearn}}

\noindent\textbf{(i)} We are given that $c=0$. For small enough $h>0$, let
$\{H_{t}\}_{t\in\lbrack0,1+h]}$ be the corresponding functions defined by
(\ref{hf}).

First prove
\begin{equation}
\lim_{n\rightarrow\infty}\left\vert E_{P^{s^{\ast}}}\left[  H_{1}\left(
\frac{\Sigma_{1}^{n}X_{i}}{\sqrt{n}}\right)  \right]  -E_{P^{\ast}}\left[
\varphi\left(  W_{1+h}^{0}\right)  \right]  \right\vert =0 \label{proof-1}%
\end{equation}
We have{\normalsize
\begin{align*}
&  E_{P^{s^{\ast}}}\left[  H_{1}\left(  \frac{\Sigma_{1}^{n}X_{i}}{\sqrt{n}%
}\right)  \right]  -E_{P^{\ast}}\left[  \varphi\left(  W_{1+h}^{0}\right)
\right] \\
=  &  E_{P^{s^{\ast}}}\left[  H_{1}\left(  \frac{\Sigma_{1}^{n}X_{i}}{\sqrt
{n}}\right)  \right]  -H_{0}(0)\\
=  &  \sum\limits_{m=1}^{n}\left\{  E_{P^{s^{\ast}}}\left[  H_{\frac{m}{n}%
}\left(  \frac{\Sigma_{1}^{m}X_{i}}{\sqrt{n}}\right)  \right]  -E_{P^{s^{\ast
}}}\left[  H_{\frac{m-1}{n}}\left(  \frac{\Sigma_{1}^{m-1}X_{i}}{\sqrt{n}%
}\right)  \right]  \right\} \\
=  &  \sum\limits_{m=1}^{n}\left\{  E_{P^{s^{\ast}}}\left[  H_{\frac{m}{n}%
}\left(  \frac{\Sigma_{1}^{m}X_{i}}{\sqrt{n}}\right)  \right]  -E_{P^{s^{\ast
}}}\left[  L_{m,n}\left(  \frac{\Sigma_{1}^{m-1}X_{i}}{\sqrt{n}}\right)
\right]  \right\} \\
&  +\sum_{m=1}^{n}\left\{  E_{P^{s^{\ast}}}\left[  L_{m,n}\left(  \frac
{\Sigma_{1}^{m-1}X_{i}}{\sqrt{n}}\right)  \right]  -E_{P^{s^{\ast}}}\left[
H_{\frac{m-1}{n}}\left(  \frac{\Sigma_{1}^{m-1}X_{i}}{\sqrt{n}}\right)
\right]  \right\} \\
=  &  :J_{1n}+J_{2n}\text{,}%
\end{align*}
}where ${L_{m,n}(x)=H_{\frac{m}{n}}(x)+\frac{\overline{\sigma}^{2}}{2n}\left[
{H_{\frac{m}{n}}^{\prime\prime}(x)}\right]  ^{+}-\frac{\underline{\sigma}^{2}%
}{2n}\left[  {H_{\frac{m}{n}}^{\prime\prime}(x)}\right]  ^{-},\ 1\leq m\leq
n}$.

By a {similar argument to that in the proof of Lemma \ref{lemma-taylor-O},
(using Lemma \ref{lemma-ddp}(3) and the fact that $E_{P^{s^{\ast}}}[X_{m}%
^{2}|\mathcal{G}_{m-1}]=I_{\{\Sigma_{1}^{m-1}X_{i}\leq0\}}\overline{\sigma
}^{2}+I_{\{\Sigma_{1}^{m-1}X_{i}>0\}}\underline{\sigma}^{2}$), deduce that
\[
\lim_{n\rightarrow\infty}|J_{1n}|=0.
\]
On the other hand, by Lemma \ref{lemma-ddp}(5), {(argue as in the proof that
$|I_{2n}|\rightarrow0$ in Appendix \ref{app-CLTproof}),} we have
$\lim_{n\rightarrow\infty}|J_{2n}|=0$. Thus we obtain (\ref{proof-1}).
%\begin{align*}
%|J_{2n}| &  \leq\sum_{m=1}^{n}E_{P^{s^{\ast}}}\left[  \left\vert L_{\frac
%{m}{n}}\left(  \frac{S_{m-1}}{\sqrt{n}}\right)  -H_{\frac{m-1}{n}}\left(
%\frac{S_{m-1}}{\sqrt{n}}\right)  \right\vert \right]  \\
%&  \leq\sum_{m=1}^{n}\sup\limits_{x\in\mathbb{R}}\left\vert L_{\frac{m}{n}%
%}(x)-H_{\frac{m-1}{n}}(x)\right\vert \\
%&  =\sum_{m=1}^{n}\sup\limits_{x\in\mathbb{R}}\left\vert H_{\frac{m-1}{n}%
%}(x)-H_{\frac{m}{n}}\left(  x\right)  -\frac{\overline{\sigma}^{2}}{2n}\left[
%{H_{\frac{m}{n}}^{\prime\prime}(x)}\right]  ^{+}+\frac{\underline{\sigma}^{2}%
%}{2n}\left[  {H_{\frac{m}{n}}^{\prime\prime}(x)}\right]  ^{-}\right\vert \\
%&  \rightarrow0\text{ \ as }n\rightarrow\infty\text{.}%
%\end{align*}
}

{By the definition of functions $\{H_{t}\}$ and Lemma \ref{lemma-ddp}(6),} and
arguing as at the end of Appendix \ref{app-CLTproof}, the proof of (i) is
complete.
%By Lemma \ref{lemma-ddp}(6) and the definition of $H_{t}$, we have
%\begin{align*}
%&  \lim_{n\rightarrow\infty}\left\vert E_{P^{s^{\ast}}}\left[  \varphi\left(
%\frac{S_{n}}{\sqrt{n}}\right)  \right]  -E_{P^{\ast}}\left[  \varphi\left(
%W_{1}^{0}\right)  \right]  \right\vert \\
%\leq &  \lim_{n\rightarrow\infty}\left\vert E_{P^{s^{\ast}}}\left[
%\varphi\left(  \frac{S_{n}}{\sqrt{n}}\right)  \right]  -E_{P^{s^{\ast}}%
%}\left[  H_{1}\left(  \frac{S_{n}}{\sqrt{n}}\right)  \right]  \right\vert \\
%&  +\lim_{n\rightarrow\infty}\left\vert E_{P^{s^{\ast}}}\left[  H_{1}\left(
%\frac{S_{n}}{\sqrt{n}}\right)  \right]  -E_{P^{\ast}}\left[  \varphi\left(
%W_{1+h}^{0}\right)  \right]  \right\vert \\
%&  +\left\vert E_{P^{\ast}}\left[  \varphi\left(  W_{1+h}^{c=0}\right)
%\right]  -E_{P^{\ast}}\left[  \varphi\left(  W_{1}^{0}\right)  \right]
%\right\vert \\
%\leq &  \sup_{x\in\mathbb{R}}\left\vert H_{1}\left(  x\right)  -\varphi\left(
%x\right)  \right\vert +C_{0}\sqrt{\overline{\sigma}^{2}+\underline{\sigma}%
%^{2}}\sqrt{h}\\
%\leq &  2C_{0}\sqrt{\overline{\sigma}^{2}+\underline{\sigma}^{2}}\sqrt{h}%
%\end{align*}
%where $C_{0}=\Vert\varphi^{\prime}\Vert$. Since $h$ is arbitrary, the proof of
%(i) is complete.

\smallskip

\noindent\textbf{(ii)} By Corollary \ref{cor-pro-O}, we have that, for any
$N>0$,{\normalsize { {
\[
P^{s^{\ast}}\left(  \cap_{n=N}^{\infty}\{\Sigma_{1}^{n}X_{i}>0\}\right)
\leq\lim_{n\rightarrow\infty}\sup_{s\in\mathcal{S}}P^{s}\left(  \Sigma_{1}%
^{n}X_{i}/\sqrt{n}>0\right)  =\frac{\overline{\sigma}}{\overline{\sigma
}+\underline{\sigma}}<1.
\]
}}}By the corresponding result for the indicator of $(-\infty,c]$, (see
Corollary 3.4 in our working paper version),{\normalsize {{
\[
P^{s^{\ast}}\left(  \cap_{n=N}^{\infty}\{\Sigma_{1}^{n}X_{i}\leq0\}\right)
\leq\lim_{n\rightarrow\infty}\sup_{s\in\mathcal{S}}P^{s}\left(  \Sigma_{1}%
^{n}X_{i}/\sqrt{n}\leq0\right)  =\frac{\overline{\sigma}}{\overline{\sigma
}+\underline{\sigma}}<1.
\]
}}}

\noindent\textbf{(iii)} To derive (\ref{relativefreq}), argue first, as in
Corollary \ref{cor-pro-O}, that the indicator for $[0,\infty)$ can be
approximated by a function $\varphi$ satisfying conditions of the CLT and the
bandit application. Then it can be shown that (\ref{s*}) is asymptotically
optimal also when the indicator replaces $\varphi$, that is, when DM solves
$\sup_{s\in\mathcal{S}}P^{s}\left(  \Sigma_{1}^{n}X_{i}/\sqrt{n}>d\right)  $.
Finally, apply the closed-form expression in the noted
corollary.{\normalsize \hfill$\blacksquare$ \medskip}

\subsection{Proof of Theorem \ref{thm-learn} and Remark
\ref{remark-posteriors}}

\noindent\textbf{Theorem }\ref{thm-learn}: Bayesian updating implies that
$\{\mu_{n}\}$ is a $P^{s^{\ast}}$-martingale adapted to $\{\mathcal{G}_{n}\}$.
Since $\{\mu_{n}\}$ is uniformly bounded, there exists a random variable $\mu$
such that{\normalsize
\[
\lim_{n\rightarrow\infty}\mu_{n}=\mu\text{ \ \ }P^{s^{\ast}}\text{-a.s.}%
\]
}

\noindent Step 1: $\mu=0$ or $1$ $P^{s^{\ast}}$-a.s., which implies
(\ref{mulimit}): Purely for simplicity, we give the argument when
$\underline{p}+\overline{p}=1$; the proof for the general case will be evident.

We have $P^{s^{\ast}}(\widehat{\Omega})=1$, where $\widehat{{\Omega}}%
=\{\omega\in\Omega\ |\ \lim_{n\rightarrow\infty}\mu_{n}(\omega)=\mu(\omega
)\}$. For any $\omega\in\widehat{{ \Omega}}$, {\normalsize
\[
\mu_{n}(\omega)=\frac{\underline{p}\mu_{n-1}(\omega)}{\underline{p}\mu
_{n-1}(\omega)+\overline{p}(1-\mu_{n-1}(\omega))}\ \text{ or }\ \frac
{\overline{p}\mu_{n-1}(\omega)}{\overline{p}\mu_{n-1}(\omega)+\underline
{p}(1-\mu_{n-1}(\omega))}.
\]
}Thus, without loss of generality, there exists a subsequence $\{\mu_{k_{n}%
}\}$ satisfying%
\[
\mu_{k_{n}}(\omega)=\frac{\underline{p}\mu_{k_{n}-1}(\omega)}{\underline{p}%
\mu_{k_{n}-1}(\omega)+\overline{p}(1-\mu_{k_{n}-1}(\omega))},
\]
which implies that{\normalsize
\[
\mu(\omega)=\frac{\underline{p}\mu(\omega)}{\underline{p}\mu(\omega
)+\overline{p}(1-\mu(\omega))}.
\]
}Thus $\mu(\omega)=0$ or $1$.

\smallskip

\noindent Step 2: For $n\geq1$, define%
\[
\underline{M}_{n}=\min\{\mu_{n},1-\mu_{n}\}\text{, }\overline{M}_{n}=\max
\{\mu_{n},1-\mu_{n}\}
\]
Then, by the dominated convergence theorem, {\normalsize
\[
\lim_{n\rightarrow\infty}E_{P^{s^{\ast}}}\left[  \underline{M}_{n}\right]
=E_{P^{s^{\ast}}}\left[  \lim_{n\rightarrow\infty}\underline{M}_{n}\right]
=E_{P^{s^{\ast}}}\left[  \min\{\mu,1-\mu\}\right]  =0.
\]
}

\medskip

{For small enough $h>0$, }let $\{H_{t}\}_{t\in\lbrack0,1+h]}$ be the functions
defined in (\ref{hf}), and let $\{L_{m,n}\}_{m=1}^{n}$ be the functions
defined in (\ref{Lt}). We prove below that
\begin{equation}
\lim_{n\rightarrow\infty}\sum_{m=1}^{n}\left\vert E_{P^{s^{\ast}}}\left[
H_{\frac{m}{n}}\left(  \frac{\Sigma_{1}^{m}X_{i}}{\sqrt{n}}\right)  \right]
-E_{P^{s^{\ast}}}\left[  L_{m,n}\left(  \frac{\Sigma_{1}^{m-1}X_{i}}{\sqrt{n}%
}\right)  \right]  \right\vert =0. \label{HL}%
\end{equation}
This is the counterpart for the present setting of the limit result
(\ref{lemma-taylor1}) in the proof of our CLT (Lemma \ref{lemma-taylor-O}),
where instead of the expectation with respect to the single measure
$P^{s^{\ast}}$, one has the upper expectation $\mathbb{E}$ corresponding to
the set of measures $\mathcal{P}$. The proof of (\ref{HL}) roughly parallels
the earlier arguments but the difference between $E_{P^{s^{\ast}}}$ and
$\mathbb{E}$ necessitates some adjustments (notably in Step 4).

Define {\normalsize
\[
d(m,n)=E_{P^{s^{\ast}}}\left[  H_{\frac{m}{n}}\left(  \frac{\Sigma_{1}%
^{m-1}X_{i}}{\sqrt{n}}\right)  +H_{\frac{m}{n}}^{\prime}\left(  \frac
{\Sigma_{1}^{m-1}X_{i}}{\sqrt{n}}\right)  \frac{X_{m}}{\sqrt{n}}+H_{\frac
{m}{n}}^{\prime\prime}\left(  \frac{\Sigma_{1}^{m-1}X_{i}}{\sqrt{n}}\right)
\frac{X_{m}^{2}}{2n}\right]  .
\]
}It suffices for (\ref{HL}) to prove that
\begin{align}
&  \sum_{m=1}^{n}\left\vert E_{P^{s^{\ast}}}\left[  H_{\frac{m}{n}}\left(
\frac{\Sigma_{1}^{m}X_{i}}{\sqrt{n}}\right)  \right]  -d(m,n)\right\vert
\rightarrow0\text{ \ and}\label{dH}\\
&  \sum_{m=1}^{n}\left\vert d(m,n)-E_{P^{s^{\ast}}}\left[  L_{m,n}\left(
\frac{\Sigma_{1}^{m-1}X_{i}}{\sqrt{n}}\right)  \right]  \right\vert
\longrightarrow0\text{.} \label{dL}%
\end{align}

\smallskip

\noindent Step 3: Prove (\ref{dH}). The argument is similar to that for
(\ref{O-var-taylor1}).

\smallskip

\noindent Step 4: Prove (\ref{dL}). By (\ref{Psn}), for any $m\geq1$,
$E_{P^{s^{\ast}}}[X_{m}|\mathcal{G}_{m-1}]=0$, and%
\begin{equation}
E_{P^{s^{\ast}}}[X_{m}^{2}|\mathcal{G}_{m-1}]=\left\{
\begin{array}
[c]{ccc}%
\overline{\sigma}^{2}\overline{M}_{m}+\underline{\sigma}^{2}\underline{M}_{m}
& \text{if} & \Sigma_{1}^{m-1}X_{i}\leq0\\
&  & \\
\underline{\sigma}^{2}\overline{M}_{m}+\overline{\sigma}^{2}\underline{M}_{m}
& \text{if} & \Sigma_{1}^{m-1}X_{i}>0
\end{array}
\right.  \label{s*variance}%
\end{equation}
Therefore, for $C_{1}$ equal to the uniform bounded of $|H_{t}^{^{\prime
\prime}}(x)|$,
\begin{align*}
&  \sum_{m=1}^{n}\left\vert d(m,n)-E_{P^{s^{\ast}}}\left[  L_{m,n}\left(
\tfrac{\Sigma_{1}^{m-1}X_{i}}{\sqrt{n}}\right)  \right]  \right\vert \\
\leq &  \sum_{m=1}^{n}E_{P^{s^{\ast}}}\left[  \frac{1}{2n}\left[  H_{\frac
{m}{n}}^{\prime\prime}\left(  \tfrac{\Sigma_{1}^{m-1}X_{i}}{\sqrt{n}}\right)
\right]  ^{+}\left(  \overline{\sigma}^{2}-\overline{\sigma}^{2}\overline
{M}_{m}-\underline{\sigma}^{2}\underline{M}_{m}\right)  \right] \\
&  +\sum_{m}^{n}E_{P^{s^{\ast}}}\left[  \frac{1}{2n}\left[  H_{\frac{m}{n}%
}^{\prime\prime}\left(  \tfrac{\Sigma_{1}^{m-1}X_{i}}{\sqrt{n}}\right)
\right]  ^{-}\left(  \underline{\sigma}^{2}\overline{M}_{m}+\overline{\sigma
}^{2}\underline{M}_{m}-\underline{\sigma}^{2}\right)  \right] \\
\leq &  \frac{C_{1}(\overline{\sigma}^{2}-\underline{\sigma}^{2})}{n}%
\sum_{m=1}^{n}E_{P^{s^{\ast}}}\left[  \underline{M}_{m}\right]  \text{
}\longrightarrow0\text{ (by Step 2).}%
\end{align*}

\begin{remark}
Step 4 involves a departure from the arguments of the CLT. In the latter, we
had by assumption (\ref{condvar}) that upper and lower conditional variances
were constant and equal to $\overline{\sigma}^{2}$ and $\underline{\sigma}%
^{2}$ respectively, while here the relevant conditional variances are under
$P^{s^{\ast}}$ and are stochastic as shown in (\ref{s*variance}). Also
noteworthy is that, while all other steps in the argument are valid for all
strategies $s$, Step 4 relies explicitly on $s=s^{\ast}$.
\end{remark}

\smallskip

\noindent Step 5: Complete the proof. It can be checked that,
\begin{align*}
&  E_{P^{s^{\ast}}}\left[  H_{1}\left(  \frac{\Sigma_{1}^{n}X_{i}}{\sqrt{n}%
}\right)  \right]  -H_{0}(0)\\
=  &  \sum\limits_{m=1}^{n}\left\{  E_{P^{s^{\ast}}}\left[  H_{\frac{m}{n}%
}\left(  \frac{\Sigma_{1}^{m}X_{i}}{\sqrt{n}}\right)  \right]  -E_{P^{s^{\ast
}}}\left[  H_{\frac{m-1}{n}}\left(  \frac{\Sigma_{1}^{m-1}X_{i}}{\sqrt{n}%
}\right)  \right]  \right\} \\
=  &  \sum\limits_{m=1}^{n}\left\{  E_{P^{s^{\ast}}}\left[  H_{\frac{m}{n}%
}\left(  \frac{\Sigma_{1}^{m}X_{i}}{\sqrt{n}}\right)  \right]  -E_{P^{s^{\ast
}}}\left[  L_{m,n}\left(  \frac{\Sigma_{1}^{m-1}X_{i}}{\sqrt{n}}\right)
\right]  \right\} \\
&  +\sum_{m=1}^{n}\left\{  E_{P^{s^{\ast}}}\left[  L_{m,n}\left(  \frac
{\Sigma_{1}^{m-1}X_{i}}{\sqrt{n}}\right)  \right]  -E_{P^{s^{\ast}}}\left[
H_{\frac{m-1}{n}}\left(  \frac{\Sigma_{1}^{m-1}X_{i}}{\sqrt{n}}\right)
\right]  \right\} \\
=  &  :\hat{J}_{1n}+\hat{J}_{2n}\text{.}%
\end{align*}
By (\ref{HL}), we have $\lim_{n\rightarrow\infty}|\hat{J}_{1n}|=0$. By Lemma
\ref{lemma-ddp}(5), (argue as in the proof that $|I_{2n}|\rightarrow0$ in
Appendix \ref{app-CLTproof}), we have $\lim_{n\rightarrow\infty}|\hat{J}%
_{2n}|=0$. Therefore,
\[
\lim_{n\rightarrow\infty}\left\vert E_{P^{s^{\ast}}}\left[  H_{1}\left(
\frac{\Sigma_{1}^{n}X_{i}}{\sqrt{n}}\right)  \right]  -H_{0}(0)\right\vert
=0.
\]
By the definition of functions $\{H_{t}\}$, with arguments similar to those at
the end of Appendix \ref{app-CLTproof}, we have
\[
\left\vert E_{P^{s^{\ast}}}\left[  \varphi\left(  \frac{\Sigma_{1}^{n}X_{i}%
}{\sqrt{n}}\right)  \right]  -E_{P^{\ast}}[\varphi(W_{1}^{0})]\right\vert
\rightarrow0\quad\text{as }\ n\rightarrow\infty\text{. \ \ \ \ }{\blacksquare}%
\]
\hfill

\bigskip

\noindent\textbf{Remark \ref{remark-posteriors}}: (i) is proven in Step 1
above. It is assumed there that $s=s^{\ast}$, but the identical arguments
apply to any $s$.

Consider (ii). Let $\nu_{n}=\mu_{n}/(1-\mu_{n})$ and apply (\ref{posteriors})
to derive, for any $s$,
\begin{align*}
&  \log\nu_{n+1}-\log\nu_{1}\\
&  =\left[  (f_{a}^{s}(\omega^{(n)})-f_{a,0}^{s}(\omega^{(n)}))-(f_{b}%
^{s}(\omega^{(n)})-f_{b,0}^{s}(\omega^{(n)}))\right]  \log\left(
\frac{\underline{p}}{\overline{p}}\right) \\
&  +[f_{a,0}^{s}(\omega^{(n)})-f_{b,0}^{s}(\omega^{(n)})]\log\left(
\frac{1-\underline{p}}{1-\overline{p}}\right)  \text{.}%
\end{align*}

{Define the sets}

{%
\begin{align*}
N_{a}=  &  \left\{  \omega:\lim_{n\rightarrow\infty}f_{a}^{s}(\omega
^{(n)})=\infty,\ \lim_{n\rightarrow\infty}f_{b}^{s}(\omega^{(n)}%
)<\infty\right\}  \text{,}\\
N_{b}=  &  \left\{  \omega:\lim_{n\rightarrow\infty}f_{a}^{s}(\omega
^{(n)})<\infty,\ \lim_{n\rightarrow\infty}f_{b}^{s}(\omega^{(n)}%
)=\infty\right\}  \text{,}\\
N_{a,b}=  &  \left\{  \omega:\lim_{n\rightarrow\infty}f_{a}^{s}(\omega
^{(n)})=\infty,\ \lim_{n\rightarrow\infty}f_{b}^{s}(\omega^{(n)}%
)=\infty\right\}  \text{,}\\
M_{a}=  &  \left\{  \omega:\lim_{n\rightarrow\infty}\frac{f_{a,0}^{s}%
(\omega^{(n)})}{f_{a}^{s}(\omega^{(n)})}=1-\underline{p}\right\}  \text{,}\\
M_{b}=  &  \left\{  \omega:\lim_{n\rightarrow\infty}\frac{f_{b,0}^{s}%
(\omega^{(n)})}{f_{b}^{s}(\omega^{(n)})}=1-\overline{p}\right\}  \text{.}%
\end{align*}
}

\noindent Consider $\omega\in N_{a,b}\cap M_{a}\cap M_{b}$: Then $\log
\nu_{n+1}-\log\nu_{1}=$
\begin{align*}
&  -f_{a}^{s}\left[  \underline{p}\log\left(  \frac{\overline{p}}%
{\underline{p}}\right)  +(1-\underline{p})\log\left(  \frac{1-\overline{p}%
}{1-\underline{p}}\right)  \right] \\
&  -f_{b}^{s}\left[  \overline{p}\log\left(  \frac{\underline{p}}{\overline
{p}}\right)  +(1-\overline{p})\log\left(  \frac{1-\underline{p}}%
{1-\overline{p}}\right)  \right] \\
&  \equiv-f_{a}^{s}H_{1}-f_{b}^{s}H_{2}\text{. }%
\end{align*}
By the concavity of $\log$, $H_{1},H_{2}<0$. Therefore, $\nu_{n}%
\rightarrow\infty$, equivalently $\mu_{n}\rightarrow1$, on $N_{a,b}\cap
M_{a}\cap M_{b}$. By the LLN, $Q^{s}(N_{a,b}\cap M_{a}\cap M_{b}%
)=Q^{s}(N_{a,b})$. Conclude that
\[
Q^{s}\left(  N_{a,b}\cap\{\omega:\mu_{n}\rightarrow1\}\right)  =Q^{s}\left(
N_{a,b}\right)  \text{.}%
\]

Similar equations apply if $N_{a,b}$ is replaced by either $N_{a}$ or $N_{b}$.
{Finally, }since $\{N_{a},N_{b},N_{a,b}\}$ is a partition of $\Omega$,
conclude that $Q^{s}(\{\omega:\mu_{n}\rightarrow1\})=1$.

\end{document}